\documentclass[preprint,12pt,times]{elsarticle}

\usepackage{amsmath}
\usepackage{amssymb}
\usepackage{amsfonts}
\usepackage{amsthm}
\usepackage{amscd}
\usepackage{graphicx}
\usepackage{xcolor} 
\usepackage{framed}
\usepackage{comment}
\usepackage{chngcntr}
\usepackage{mathrsfs}
\usepackage{titletoc}
\usepackage{bm}
\usepackage[pagewise,mathlines]{lineno}

\biboptions{sort&compress}

\newtheorem{theorem}{Theorem}
\newtheorem{proposition}[theorem]{Proposition}
\newtheorem{lemma}[theorem]{Lemma}

\newtheorem{remark}[theorem]{Remark}

\makeatletter

\newcommand{\Rmnum}[1]{\expandafter\@slowromancap\romannumeral #1@}
\makeatother

\newif\ifrevisionmarks

\revisionmarkstrue        

\ifrevisionmarks

\else

\fi

\usepackage[breaklinks]{hyperref}

\ifrevisionmarks
  \hypersetup{
    colorlinks=true,
    allcolors=blue
  }
\else
  \hypersetup{hidelinks}
\fi

\journal{***}

\begin{document}

\title{Fractional Sobolev spaces and fractional $p$-Laplace equations on locally finite graphs}

\author[author1]{Yong Lin}
\ead{yonglin@mail.tsinghua.edu.cn}

\author[author2]{Yunyan Yang}
\ead{yunyanyang@ruc.edu.cn}

\author[author3]{Mengjie Zhang\corref{cor1}}
\ead{zhangmengjie@sdufe.edu.cn}

\cortext[cor1]{Corresponding author.}

\address[author1]{Yau Mathematical Sciences Center,
Tsinghua University, Beijing 100084, China}

\address[author2]{School of Mathematics,
Renmin University of China, Beijing 100872, China}

\address[author3]{Key Laboratory of Blockchain Finance,
Shandong University of Finance and Economics,\\
Jinan 250014, China}

\begin{frontmatter}

\begin{abstract}

In this paper, we develop a theory of fractional Sobolev spaces  for $s\in(0,1)$ and $p\in[2,\infty)$
on locally finite, connected and stochastically complete weighted graphs with infinitely many vertices.
Assuming a uniform positive lower bound on the vertex measures, we establish completeness, reflexivity,
continuous embeddings and density of finitely supported functions.
Furthermore, we extend the fractional $p$-Laplace operator
from finite graphs to the present setting of infinite graphs, and establish an integration by parts formula.
Under suitable assumptions on the potential and nonlinearity,
we establish the existence of strictly positive solutions and ground state solutions to the associated nonlinear Schr\"odinger equations.
 These results provide an analytical framework for the variational study of nonlinear fractional equations on locally finite graphs.

\end{abstract}



\end{frontmatter}

\section{Introduction}

Fractional Laplace operators arise naturally in potential theory,
fractional calculus and probability, and have important applications
to differential equations and mathematical physics
\cite{28,29,30,31}.
Their analysis is closely connected with fractional Sobolev spaces,
which provide natural settings for many nonlocal equations.
For background on Sobolev spaces, fractional operators and related
problems, we refer to \cite{A1,A,A78,A12,Kwa15,del-2021}
and the references therein.

On lattices, the linear fractional Laplacian and related equations
have been studied from several perspectives.
These include connections with superdiffusion
\cite{CRSTV2015} and one-dimensional operators \cite{LM22},
regularity and almost periodicity for semidiscrete evolution
equations \cite{LR2018}, and optimal Hardy inequalities
\cite{KN23}.
Ciaurri et al.\ \cite{CRSTV2018} obtained kernel estimates
and regularity results, and studied the approximation of
continuous fractional problems by discrete ones.
They also established fractional discrete Sobolev and
Poincar\'e inequalities and developed a Hilbert space
formulation for nonlocal Dirichlet problems.

A general account of analysis on discrete weighted graphs
can be found in the monograph \cite{Keller-1}.
Shao, Yang and Zhao \cite{S-Y-Z-1} developed a general theory
of Sobolev spaces on locally finite graphs, including
completeness, reflexivity, separability and Sobolev embeddings.
Their approach uses vector-valued function spaces whose dimensions vary from vertex to vertex to accommodate the local graph gradient.

 The density of finitely supported functions in first-order
Sobolev spaces has also been studied on graphs.
Under a uniform positive lower bound on the vertex measures
and a uniform upper bound on the weighted degrees
$\sum_{y\sim x}w_{xy}$,
Han, Shao and Zhao \cite{Han-Shao-Zhao2020} proved that
$W^{1,2}_0(V)=W^{1,2}(V)$.
Under the same bounds, Shao, Yang and Zhao
\cite{S-Y-Z-1} established
$W^{1,p}_0(V)=W^{1,p}(V)$ for $1\leq p<\infty$.
Tian and Zhao \cite{Tian-Zhao2025} showed that this equality
can fail under different assumptions on the graph and
the vertex measure.
These results motivate the study of density in the
fractional setting considered here.

Nonlinear equations on graphs have also been studied by using
variational and topological methods.
Grigor'yan, Lin and Yang
\cite{A-Y-Y-1,A-Y-Y-2,A-Y-Y-3} established existence results for several classes of nonlinear equations by variational methods.
Further developments concern nonlinear Schr\"odinger equations
\cite{Hua-Xu,Zhang-Zhao}, $p$-Laplacian equations
\cite{Han-Shao}, biharmonic Choquard equations
\cite{Liu-Zhang}, Kazdan--Warner equations \cite{Ge,S-w},
and Chern--Simons--Higgs equations \cite{HS}.
In the setting of metric graphs, De Coster et al.\ \cite{De-C}
 compared ground states defined by minimization on the Nehari manifold with least action solutions.

 Building on the $p=2$ theory developed in
\cite{Z-L-Y2}, our first aim is to establish
a fractional Sobolev framework for $p\in [2,\infty)$
on locally finite, connected and stochastically complete weighted graphs with infinitely many vertices. We assume that the vertex measures have a uniform
positive lower bound.
For $s\in(0,1)$, we establish completeness, reflexivity
and continuous embeddings for $W^{s,p}(V)$.
 We also prove that $C_c(V)$ is a proper dense subspace of
$W^{s,p}(V)$, and hence $W_0^{s,p}(V)=W^{s,p}(V)$.
This density result holds under the stated assumptions,
without requiring a uniform upper bound on
the weighted degrees.
The fractional gradient involves interactions with infinitely many vertices.
We handle this feature by combining fractional kernel summability estimates with vector-valued function spaces.

Our second aim is to extend the fractional $p$-Laplace
operator studied on finite graphs in \cite{Yu2026}
to the present infinite graph setting.
This also generalizes the fractional Laplace operator
considered in \cite{Z-L-Y2} from $p=2$ to
$p\in[2,\infty)$.
The operator is associated with the energy
$  p^{-1}\int_V|\nabla^s u|^p d\mu$
and is nonlinear when $p>2$.
We show that the operator is well defined and establish
an integration by parts formula.
Under suitable assumptions on the potential and nonlinearity,
we use this variational framework to obtain
strictly positive solutions and ground state solutions
of the associated nonlinear Schr\"odinger equations.

\section{Notations and main results}

\subsection{Notations}

Let $G=(V,E,\mu,w)$ be a locally finite, connected weighted
graph with the infinite  vertex set $V$ and the edge set $E$.
The graph is undirected, with a positive vertex measure $\mu:V\rightarrow \mathbb{R}^+$ and a positive symmetric edge weight function $w:E\rightarrow \mathbb{R}^+$.

Denote by $C(V)$ the space of all real-valued functions
on $V$ and by $C_c(V)$ its subspace of finitely supported
functions.  Let   $\mathrm{supp} (u)=\{x\in V: u(x)\not=0\}$ be the support of   $u$.
Integration with respect to $\mu$ is understood as $\int_V u\,d\mu=\sum_{x\in V}u(x)\mu(x)$.
 For $q\in[1,\infty)$, we define
$$L^q(V)=\left\{u\in C(V):\int_V|u|^q\,d\mu<\infty \right\}$$
and equip it with the norm
$\|u\|_{L^q(V)}= (\int_V|u|^q\,d\mu)^{1/q}$.
We also write $L^\infty(V)$ for the space of all bounded
real-valued functions on $V$, equipped with the norm
$\|u\|_{L^\infty(V)}=\sup_{x\in V}|u(x)|$.
We write $$\langle u,v\rangle =\int_V uv\,d\mu$$ whenever
$uv\in L^1(V)$, and use the same notation for duality pairings.

 The formal  graph Laplacian is defined by
\begin{align*}
\Delta u(x)
=\frac{1}{\mu(x)}
\sum_{y\sim x}w_{xy}\bigl(u(y)-u(x)\bigr),\quad \forall u\in C(V),
\end{align*}
where $y\sim x$ means that $y$ and $x$ are adjacent.
Following
\cite[Chapter~1, Section~1]{Keller-1}, let $\Delta_D$ denote
the nonpositive self-adjoint Dirichlet realization
of $\Delta$ on $L^2(V)$, and let
$P_t=\mathrm{e}^{t\Delta_D}$ be the associated heat semigroup.
As explained in \cite[Chapter~5, Section~1]{Keller-1},
this semigroup admits a heat kernel with respect
to $\mu$, given by
\[
p(t,x,y)=\frac{(P_t\delta_y)(x)}{\mu(y)},
\quad \forall t\ge0,\  x,y\in V,
\]
where $\delta_y$ is the indicator function of the
vertex $y$.
The kernel is nonnegative and symmetric in $x$ and $y$,
and satisfies  $\sum_{z\in V}p(t,x,z)\mu(z)\leq1$ for every
$ t>0$ and $x \in V$.
Applying \cite[Theorem~A.33]{Keller-1} to the
nonnegative self-adjoint operator $-\Delta_D$,
and using the continuity of point evaluation
on $L^2(V)$, we obtain the heat equation
\begin{align}\label{heat-equation}
\frac{\partial}{\partial t}p(t,x,y)
=
\Delta_x p(t,x,y)
=
\Delta_y p(t,x,y),
\quad \forall t>0,\ x,y\in V,
\end{align}
where $\Delta_x$ and $\Delta_y$ denote the formal
graph Laplacian acting on the first and second
spatial variables, respectively.
 The graph $G$ is called stochastically complete if
\begin{align}\label{comp}
\sum_{y\in V}p(t,x,y)\mu(y)=1,
\quad \forall t>0,\  x\in V.
\end{align}
Otherwise, $G$ is called stochastically incomplete.
See also
\cite{Haeseler-Keller-Lenz-Wojciechowski,Keller-Lenz,Huang,Hua-Lin2017}.
Throughout this paper, we assume that $G$ is
stochastically complete.

 Since $G$ is locally finite and connected and has
infinitely many vertices, its vertex set $V$ is countably
infinite. We fix an enumeration
$$
V=\{x_1,x_2,\cdots,x_i,\cdots\}.
$$
Denote by $\mathbb R^\infty$ the vector space of all
real sequences, namely,
$$
\mathbb R^\infty
=\left\{\mathbf a=(a_1,a_2,\cdots, a_i, \cdots):
a_i\in\mathbb R\ \text{ for every }\ i\ge1\right\}.
$$
The space of all $\mathbb R^\infty$-valued functions
on $V$ is denoted by
\begin{align*}
C(V;\mathbb R^\infty)
=\left\{\mathbf{u} =(u_1 , u_2, \cdots, u_i , \cdots):
 u_i\in C(V)\ \text{ for every }\ i\ge1\right\}.
\end{align*}
Define the vector-valued function space
 \begin{align}\label{e-L}
 \mathcal{L}^q(V)=\left\{ \mathbf{u} \in  C(V;\mathbb R^\infty) : \int_V|\mathbf{u}(x)|^q d \mu<\infty\right\}, \quad \forall   q\in [1,\infty),
 \end{align}
where
$|\mathbf{u}(x)|= (\sum_{i=1}^{\infty}|u_i(x)|^2 )^{1/2}$
 for $\mathbf{u}(x)=(u_1(x),u_2(x),\cdots, u_i(x),\cdots)$. We equip $\mathcal{L}^q(V)$ with the norm
$$
\|\mathbf{u}\|_{\mathcal{L}^q(V)}
=
\left(\int_V|\mathbf{u}(x)|^q\,d\mu\right)^{\frac{1}{q}}.
$$
Let $\mathbf{e}_i\in\mathbb{R}^{\infty}$ denote the vector
whose $i$-th component is $1$ and whose other components
are all $0$.
For any fixed vertex $x_0\in V$,
the vector-valued function
$$
\boldsymbol{\varphi}_i(x)
=
\begin{cases}
\mathbf{e}_i, & x=x_0,\\
\mathbf{0}, & x\neq x_0,
\end{cases}
$$
belongs to $\mathcal{L}^q(V)$ and satisfies
$\|\boldsymbol{\varphi}_i\|_{\mathcal{L}^q(V)}^q
=\mu(x_0).$
In particular, $\mathcal{L}^q(V)$ contains nonzero
elements.
With the above norm, $\mathcal{L}^q(V)$ is a normed
linear space.
Moreover,
Lemma~\ref{L3.2} below shows that it is a Banach space
for $q\in [1,\infty)$ and is reflexive for $q\in (1,\infty)$.

\subsection{Fractional gradient operator}

Assume $s\in(0,1)$.  Following the heat kernel construction in \cite{Z-L-Y1, Z-L-Y2}, we define the fractional kernel by
\begin{align}\label{Hs-1}
W_s(x, y)= \frac{s} {\Gamma(1-s)} \mu(x)\mu(y)\int_0^{\infty}  p(t, x, y) t^{-1-s} d t,\quad \forall x \neq y \in  V,
\end{align}
 where $ \Gamma$  denotes the Gamma function. For related constructions on lattice graphs, see
\cite{Wang-frac}.  The symmetry and strict positivity of the heat kernel imply
that
\begin{align}\label{symmetry}
 W_s(x,y)=W_s(y,x)>0, \quad \forall x \neq y \in  V.
\end{align}
 Moreover, from Lemma \ref{L1} below, we have
\begin{align}\label{kernel}
 W_s(x,y)\le
\sum_{z\in V\setminus\{x\}}W_s(x,z)<\infty, \quad \forall x \neq y \in  V.
\end{align}
For convenience, we define $\widetilde W_s:V\times V\rightarrow[0,\infty)$ by
\begin{align*}
\widetilde W_s(x,y)
=
\begin{cases}
W_s(x,y), & x\ne y,\\
0, & x=y.
\end{cases}
\end{align*}

The fractional gradient operator
$\nabla^s:L^\infty(V)\rightarrow C(V;\mathbb R^\infty)$
is defined by
\begin{align}\label{e-7}
\nabla^s u(x)=\left(\sqrt{\frac{\widetilde{W}_s\left(x, x_1\right)}{2 \mu(x)}}\left(u(x)-u\left(x_1\right)\right),\cdots,\sqrt{\frac{\widetilde{W}_s\left(x, x_i\right)}{2 \mu(x)}}\left(u(x)-u\left(x_i\right)\right), \cdots\right).
\end{align}
From \eqref{kernel} and the boundedness of $u$,
 $\nabla^s u(x)$ has square-summable components
for every $x\in V$.
For $u,v\in L^\infty(V)$, we denote the inner product
of $\nabla^s u(x)$ and $\nabla^s v(x)$ by
\begin{align}\label{gra}
 \nabla^su \nabla^sv   (x)
   =   \frac{1}{2\mu(x)}\sum_{y \in V \backslash \{x\}} W_s(x, y) (u(x)-u(y))(v(x)-v(y)),
  \end{align}
  and the norm of $\nabla^s u(x)$ is given by
 \begin{align}\label{N}
|\nabla^su|  (x) =   \left(\nabla^s u  \nabla^s u(x)\right)^{\frac{1}{2}} =\left( \frac{1}{2\mu(x)}\sum_{y \in V \backslash \{x\}} W_s(x, y) (u(x)-u(y))^2 \right)^{\frac{1}{2}}.
  \end{align}
The square summability of the gradient components and
the Cauchy--Schwarz inequality ensure that both series
converge absolutely.
We next introduce the corresponding fractional Sobolev spaces.

\subsection{Fractional Sobolev space}

Throughout the remainder of this paper, we assume that
\begin{align}\label{mu}
s\in (0,1),\quad p\in [2,\infty),\quad \mu_0 =\inf_{x\in V}\mu(x)>0.
\end{align}
We define the fractional Sobolev space
\[
W^{s,p}(V)
=
\left\{
u\in L^\infty(V):
\int_V \bigl(|\nabla^s u|^p+|u|^p\bigr)\,d\mu<\infty
\right\},
\]
equipped with the norm
\[
\|u\|_{W^{s,p}(V)}
=
\left(
\int_V \bigl(|\nabla^s u|^p+|u|^p\bigr)\,d\mu
\right)^{1/p},
\]
where $|\nabla^s u|$ is defined in \eqref{N}.
Let $W^{s,p}_0(V)$ denote the completion of $C_c(V)$
with respect to this norm.
The following theorem summarizes the basic properties
of $W^{s,p}(V)$.

\begin{theorem}\label{Sobolev1}
Let $G=(V,E,\mu,w)$ be a locally finite, connected and
stochastically complete weighted graph with infinitely
many vertices.
Assume that \eqref{mu} holds.
Then the following statements hold:  \\ [0.75ex]
$\mathrm{(i)}$ $W^{s,p}(V)$  is a reflexive  Banach space;\\ [0.75ex]
$\mathrm{(ii)}$  $C_c(V)$ is a proper dense subspace of $W^{s,p}(V)$,
and hence $W^{s,p}_0(V)=W^{s,p}(V)$;\\[0.75ex]
$\mathrm{(iii)}$  $W^{s,p}(V)$ is  continuously embedded in $L^q(V)$ for all $ q\in [p,\infty]$;\\ [0.75ex]
$\mathrm{(iv)}$ for every $u\in W^{s,p}(V)$,  the functions
$u^+$, $u^-$ and $|u|$ belong to $W^{s,p}(V)$, where $u^+=\max\{u,0\}$, $u^-=\min\{u,0\}$
 and $|u|=u^+-u^-$.
\end{theorem}

\begin{remark}
 Throughout this paper, we assume $p\in[2,\infty)$.
This restriction guarantees
$C_c(V)\subseteq W^{s,p}(V)$ under the present assumptions.
If the same definition of $W^{s,{q}}(V)$ is used for
$1<{q}<2$, this inclusion may fail.

For example, consider the standard lattice $V=\mathbb Z$
with unit vertex measures and unit edge weights.
This graph satisfies all the standing graph assumptions.
By \cite[Theorem~1.1(a),(b)]{CRSTV2018},
with the lattice spacing set to one, there exist constants
$0<c_s\leq C_s<\infty$, depending only on $s$, such that
$$
\frac{c_s}{|n|^{1+2s}}
\leq W_s(n,0)=K_s(n)
\leq \frac{C_s}{|n|^{1+2s}},
\quad \forall n\in\mathbb{Z}\setminus\{0\}.
$$
Let $u(0)=1$ and $u(n)=0$ for $n\neq0$.
Then $u\in C_c(\mathbb{Z})$, $\|u\|_{L^{q}}^{q}=1$, and
$$
|\nabla^s u(n)|^{q}
=
2^{-\frac{q}{2}}K_s(n)^{\frac{q}{2}},
\quad \forall n\in\mathbb{Z}\setminus\{0\}.
$$
Consequently, if ${q}(1+2s)\leq2$, then
$$
\sum_{n\in\mathbb{Z}\setminus\{0\}}|\nabla^s u(n)|^{q}
 =
2^{-\frac{q}{2}}\sum_{n\in\mathbb{Z}\setminus\{0\}}K_s(n)^{\frac{q}{2}}
 \geq
2^{-\frac{q}{2}}c_s^{\frac{q}{2}}
\sum_{n\in\mathbb{Z}\setminus\{0\}}|n|^{-\frac{q}{2}(1+2s)}
=\infty.
$$
Hence $u\in C_c(\mathbb{Z})$ but
$u\notin W^{s,{q}}(\mathbb{Z})$ whenever ${q}(1+2s)\leq2$.
In particular, one may take $s=1/4$ and ${q}=4/3$.
\end{remark}

\subsection{Fractional $p$-Laplace operator}

For $s\in(0,1)$ and $p\in[2,\infty)$, we define
$(-\Delta)^s_p:W^{s,p}(V)\to C(V)$ by
\begin{align}\label{L-p}
(-\Delta)^s_pu(x)
=
\frac{1}{2\mu(x)}
\sum_{y\in V\setminus\{x\}}
W_s(x,y)
\left(
|\nabla^s u|^{p-2}(x)
+
|\nabla^s u|^{p-2}(y)
\right)
\left(u(x)-u(y)\right),
\end{align}
where $W_s$ and $|\nabla^s u|$ are defined
in \eqref{Hs-1} and \eqref{N}, respectively.
 For $p=2$, we adopt the convention
$|\nabla^s u|^{p-2}(x)=1$ for every $x\in V$.
With this convention, $(-\Delta)^s_p$ reduces to the
fractional Laplace operator  defined in \cite{Z-L-Y2}, namely,
\begin{align*}
(-\Delta)^s u(x)
=
\frac{1}{\mu(x)}
\sum_{y\in V\setminus\{x\}}
W_s(x,y)\left(u(x)-u(y)\right).
\end{align*}

 To verify that \eqref{L-p} is well defined, let
$u\in W^{s,p}(V)$.
From Theorem \ref{Sobolev1}(iii) and $\mu(x)\ge\mu_0>0$, we have $u\in L^\infty(V)$ and
$$
\sup_{x\in V}|\nabla^s u|(x)
\leq
\mu_0^{-\frac{1}{p}}
\|\nabla^s u\|_{\mathcal L^p}
<\infty.
$$
Hence both $u$ and $|\nabla^s u|^{p-2}$ are bounded.
Together with \eqref{kernel}, these bounds imply that
the series in \eqref{L-p} converges absolutely at every vertex.

 Moreover, the operator $(-\Delta)^s_p$ is the variational operator
associated with the energy $\int_V|\nabla^s u|^p\,d\mu/p$.
For $p>2$, this energy generally differs from the pairwise difference energy
\begin{align*}
\sum_{x\in V}\sum_{y\in V\setminus\{x\}}
W_s(x,y)|u(x)-u(y)|^p.
\end{align*}
In our construction, the term ``fractional'' refers
to the heat kernel definition of $\nabla^s$, rather than
to taking a fractional power of the nonlinear graph
$p$-Laplacian.

We next establish an integration by parts formula for use in variational arguments.
For this purpose,   we adopt a duality formulation of the fractional divergence introduced in  \cite{Z-L-Y2}.

 Fix $p\in[2,\infty)$ and write $p'=p/(p-1)$.
    For   any  $\mathbf{u} \in    \mathcal{L}^{p'}(V) $, we define
$\mathrm{div}_s \mathbf{u}\in(W^{s,p}(V))^*$ by
\begin{align}\label{div}
 \left\langle\mathrm{div}_s \mathbf{u},\varphi\right\rangle
  =-\int_V\mathbf u\cdot\nabla^s\varphi\,d\mu,
 \quad \forall \varphi\in W^{s,p}(V).
\end{align}
This is well defined, since H\"older's inequality gives
$$
 \left|\left\langle\operatorname{div}_s\mathbf u,\varphi\right\rangle\right|
 \leq \|\mathbf u\|_{\mathcal L^{p'}}
       \|\nabla^s\varphi\|_{\mathcal L^p}
 \leq \|\mathbf u\|_{\mathcal L^{p'}}\|\varphi\|_{W^{s,p}}.
$$
By taking
$\varphi=\delta_x/{\mu(x)}\in C_c(V)\subseteq W^{s,p}(V)$ in \eqref{div}, we obtain the pointwise representative of
$\mathrm{div}_s\mathbf u$ as
\begin{align*}
\mathrm{div}_s  \mathbf{u} (x)=-\frac{1}{\mu(x)}
\sum_{y\in V}\mu(y)\mathbf u(y)\cdot\nabla^s\delta_x(y),\quad \forall  \mathbf{u}  \in   \mathcal{L}^{p'}(V)  .
\end{align*}
The series is absolutely convergent by H\"older's inequality.

 For $u\in W^{s,p}(V)$, we have
\begin{align*}
\left\||\nabla^s u|^{p-2}\nabla^s u\right\|_{\mathcal L^{p'} }^{p'}
=
\left\|\nabla^s u\right\|_{\mathcal L^p}^p
<\infty,
\end{align*}
and then  $|\nabla^s u|^{p-2}\nabla^s u\in\mathcal L^{p'}(V)$.
A direct calculation yields the pointwise identity
\begin{align*}
 (-\Delta )_p^s u(x) =-\mathrm{div}_s\left( |\nabla ^su|^{p-2}\nabla ^su \right)(x)  ,\quad \forall x\in  V.
\end{align*}
Combining this identity with \eqref{div}, we obtain the following integration by parts formula.

\begin{proposition}
For every $u\in W^{s,p}(V)$ and $\varphi\in C_c(V)$,
we have
\begin{align} \label{part}
\int_{ V} \varphi (-\Delta)_p^s u  d\mu=  \int_{ V}    |\nabla^s u |^{p-2}\nabla^s u \nabla^s \varphi  d\mu .
\end{align}
\end{proposition}

At this point, we have established a theoretical framework for studying discrete fractional $p$-Laplace equations.  We next turn to analyzing  the existence of solutions to some nonlinear Schr\"odinger equations involving this nonlocal operator.

\subsection{Fractional $p$-Laplace equation}

We consider the  fractional $p$-Laplace equation
\begin{align}\label{eS1}
(-\Delta)_p^su + h|u|^{p-2}u = f(x, u) \ \ \text{ in }\   V,
\end{align}
 where $ h: V \rightarrow \mathbb{R} $ satisfies
\begin{align}\label{H-1} h_0=\inf_{x \in V}h(x)>0
\end{align}
 and for any fixed $x_0\in V$,
\begin{align}\label{H-2}
 h (x)\rightarrow \infty \  \text{ as }\ \mathrm{d}(x,x_0) \rightarrow \infty .
\end{align}
Here $\mathrm{d}(x,x_0)$   represents the minimum number of edges required for the path between
 $x$ and $x_0$.
Let  $\mathcal{H}_{s,p}(V)$  be the  completion of $C_c(V)$ under the norm
 \begin{align}\label{H}
\|u\|_{\mathcal{H}_{s,p}(V)}=\left(\int_{V} \left(|\nabla^s u|^p+h |u|^p\right ) d \mu\right)^{\frac{1}{p}}.
\end{align}
Set
$$\lambda_p =\inf _{u\in\mathcal{H}_{s,p},u\not\equiv 0} \frac{\int_{V} \left(|\nabla^s u|^p+h |u|^p\right ) d \mu}{\int_V |u|^pd\mu }.$$
Testing the defining quotient with any single vertex
function gives $$0<h_0\leq\lambda_p<\infty.$$
Moreover,  we impose the following assumptions      on $ f: V \times \mathbb{R} \rightarrow \mathbb{R} $:  \\ [0.75ex]
$\left(\mathrm{A}_1\right)$
for any $x\in V$, the function $y\mapsto f(x,y)$ is continuous on $[0,\infty)$, and $f(x,0)=0$;\\  [0.75ex]
$\left(\mathrm{A}_2\right)$   for any   constant $M>0$, there exists a constant $C_M$ such that
$ f(x, y)  \leq  C_M$ for any $  (x,y) \in V\times[0, M]$; \\  [0.75ex]
$\left(\mathrm{A}_3\right)$ there exists a constant  $\alpha>p$  such that
$0 <\alpha F(x, y)  \leq y f(x, y)$ for any $ (x,y) \in V\times \mathbb{R}^+  ,$
where $F(x, y)= \int_0^y f(x, t) d t$  is the primitive function of $f$;\\ [0.75ex]
$(\mathrm{A}_4)$   it holds that
\begin{align*} \limsup_{y \rightarrow 0^+} \frac{ f(x, y)}{y^{p-1}}<\lambda_p
\end{align*}
uniformly in $x\in V$;\\ [0.75ex]
$(\mathrm{A}_5)$ for any $x\in V$,
 $f(x,y)=0$ for all $y< 0$, and
the function $y\mapsto f(x,y)/y^{p-1}$ is strictly increasing
on $(0,\infty)$.

 Concerning the existence of positive solutions or ground state solutions to \eqref{eS1},
we have the following two theorems:

\begin{theorem}\label{T3}
Let $G=(V,E,\mu,w)$ be a locally finite, connected and stochastically complete weighted graph with the infinite
vertex set $V$.
Assume that $s\in (0,1)$, $p\in [2,\infty)$ and $\inf_{x \in V}\mu(x)>0$,
 that $ h: V \rightarrow \mathbb{R} $   satisfies \eqref{H-1} and \eqref{H-2},
 and that $ f: V \times \mathbb{R} \rightarrow \mathbb{R} $ satisfies  $(\mathrm{A}_1)$--$(\mathrm{A}_4)$.
Then the  fractional $p$-Laplace equation \eqref{eS1} admits a strictly positive solution in   ${\mathcal{H}_{s,p}(V)}$.
\end{theorem}

\begin{theorem}\label{T4}
Under the assumptions of Theorem \ref{T3}, if $f$ also satisfies
$(\mathrm{A}_5)$,
then \eqref{eS1} has a strictly positive ground state solution
in $\mathcal{H}_{s,p}(V)$.
\end{theorem}

If   $u\in\mathcal{H}_{s,p}(V)$ satisfies
	\begin{align}\label{weak-s}
		\int_{V}\left( |\nabla^s u|^{p-2}\nabla^s u\nabla^s\varphi+h|u|^{p-2}u\varphi\right) d\mu=\int_{V}f(x,u)\varphi d\mu
	\end{align}
for any $ \varphi\in  \mathcal{H}_{s,p}(V)$, then $u $ is a weak solution    of \eqref{eS1}.
Taking $\varphi=\delta_x$ in the preceding identity
and using the integration by parts formula \eqref{part},
we see that $u$ is also a pointwise solution.
From Lemma  \ref{L4.2}  below, within $\mathcal H_{s,p}(V)$, the strictly positive solutions of  \eqref{eS1}
are equivalent to the nontrivial weak solutions of
  \begin{align}\label{eS1-1}
  (-\Delta)_p^su + h|u|^{p-2}u = f(x,u^+) \ \ \text{ in }\   V.
\end{align}

 Define the energy functional $E_{s,p} :\mathcal{H}_{s,p}(V)\rightarrow\mathbb{R}$ associated with  \eqref{eS1-1} as
	\begin{align}\label{E}
		E_{s,p} (u)=\frac{1}{p}\int_{V}\left(|\nabla^s u|^p+h|u|^p\right)d\mu-\int_{V}F(x,u^+)d\mu.
	\end{align}
 Under assumptions $(\mathrm{A}_1)$--$(\mathrm{A}_4)$,
Lemma~\ref{Le4-1} below
and a standard differentiation argument show that
$E_{s,p}$ is well defined and belongs to
$C^1(\mathcal H_{s,p}(V);\mathbb R)$, with the derivative given by
\begin{equation}\label{E-C1-derivative}
\langle E_{s,p}'(u),\varphi\rangle
=
\int_V\Bigl(
|\nabla^s u|^{p-2}\nabla^s u\nabla^s\varphi
+h|u|^{p-2}u\varphi
-f(x,u^+)\varphi
\Bigr)\,d\mu
\end{equation}
for any $u,\varphi\in\mathcal H_{s,p}(V)$.
In particular, the critical points of $E_{s,p}$ are precisely
the weak solutions of \eqref{eS1-1}.
We therefore seek nontrivial critical points of $E_{s,p}$
in the following discussion.

 We prove Theorems \ref{T3} and \ref{T4} using the
mountain-pass theorem and the Nehari manifold method,
respectively.
Compactness is essential in both arguments.
For any $q\in[1,\infty)$, the uniform positive lower bound
on the vertex measures ensures that the $L^q$ norm
controls the $L^\infty$ norm.
No analogous estimate holds for general $L^q$ functions
on $\mathbb{R}^d$ with Lebesgue measure.
However, the continuous embedding of $W^{s,p}(V)$ into $L^p(V)$ need not be compact on an infinite graph.
For example, on a standard lattice with unit edge weights
and unit vertex measures, Dirac functions supported at
vertices escaping to infinity form a bounded sequence
in $W^{s,p}(V)$ with no strongly convergent subsequence
in $L^p(V)$.
This loss of compactness through escape to infinity
also occurs in continuous problems on unbounded domains.

 For the space $\mathcal H_{s,p}(V)$ considered here,
the assumption that $h$ tends to infinity provides
uniform control of the $L^p$ tails of bounded sequences.
Combining this tail control with convergence on finite
graph distance balls shows that every bounded sequence
has a subsequence that converges strongly in $L^p(V)$.
The uniform positive lower bound on the vertex measures
then upgrades this convergence to strong $L^q$ convergence
for every $q\in[p,\infty]$.
The resulting compact embeddings, stated in
Lemma~\ref{Le4-1}, provide the strong convergence
needed to verify the Palais--Smale condition
and to establish the existence of a minimizer
of the energy on the Nehari manifold.

The remainder of this paper is organized as follows.
 In Section \ref{S3},  we  prove  Theorem \ref{Sobolev1}.
 In Section  \ref{S4},
 we prove  Theorems \ref{T3}  and  \ref{T4}  respectively.
For brevity, subsequences are relabelled when no confusion can arise,
and we omit the underlying domains from the notation for
function spaces when they are clear from the context.

\section{Fractional Sobolev spaces}\label{S3}

In this section, we investigate the properties of the fractional Sobolev space $W^{s,p}(V)$.
 We first establish a summability estimate for the fractional kernel $W_s$.

\begin{lemma}\label{L1}
 We have
\begin{align*}
 W_s(x,y)\le
\sum_{z\in V\setminus\{x\}}W_s(x,z)<\infty, \quad \forall x \neq y \in  V.
\end{align*}
\end{lemma}

\begin{proof}
Fix $x\in V$ and set
\begin{align*}
d_w(x) =\frac{1}{\mu(x)}\sum_{y\sim x}w_{xy}<\infty.
\end{align*}
The heat equation  \eqref{heat-equation}  and the nonnegativity of the heat kernel imply that
\begin{align*}
\frac{\partial}{\partial t}p(t,x,x)
&=
\frac{1}{\mu(x)}
\sum_{y\sim x}w_{xy}\bigl(p(t,y,x)-p(t,x,x)\bigr)\\
&\ge -d_w(x)p(t,x,x)
\end{align*}
 for every $ t>0$.
Thus, for every fixed $x\in V$, the function $t\mapsto \mathrm{e}^{d_w(x)t}p(t,x,x)$ is nondecreasing on $(0,\infty)$.
Since $\mu(x)p(t,x,x)\to1$ as $t\to0^+$, integration gives
\begin{align*}
\mu(x)p(t,x,x)\ge \mathrm{e}^{-d_w(x)t}.
\end{align*}
Together with $0\le\mu(x)p(t,x,x)\le1$, this yields
\begin{align}\label{0-1}
0\le1-\mu(x)p(t,x,x)
\le1-\mathrm{e}^{-d_w(x)t}
\le d_w(x)t,
\quad \forall t>0.
\end{align}
Using \eqref{Hs-1}, Tonelli's theorem and stochastic completeness   (\ref{comp}), we obtain
\begin{align}\label{kernel-1}
\sum_{y\in V\setminus\{x\}}W_s(x,y)
&=\frac{s\mu(x)}{\Gamma(1-s)}
\sum_{y\in V\setminus\{x\}}\mu(y)
\int_0^\infty p(t,x,y)t^{-1-s}\,dt
\nonumber\\
&=\frac{s\mu(x)}{\Gamma(1-s)}
\int_0^\infty
\left(\sum_{y\in V\setminus\{x\}}
p(t,x,y)\mu(y)\right)t^{-1-s}\,dt
\nonumber\\
&\le\frac{s\mu(x)}{\Gamma(1-s)}
\left(
\int_0^1\bigl(1-p(t,x,x)\mu(x)\bigr)t^{-1-s}\,dt
+\int_1^\infty t^{-1-s}\,dt
\right)
\nonumber\\
&\le\frac{\mu(x)}{\Gamma(1-s)}
\left(\frac{s}{1-s} d_w(x) +1\right)
<\infty.
\end{align}
Here the interchange of the sum and the integral is
justified by the nonnegativity of the heat kernel,
and the last estimate follows from \eqref{0-1}
and $s\in(0,1)$.
  In particular, \eqref{kernel-1} implies that $W_s(x,y)$
is finite for all distinct $x,y\in V$.
\end{proof}

To prove Theorem~\ref{Sobolev1}(i), we first establish
the completeness and reflexivity of $\mathcal{L}^q $ defined by \eqref{e-L}.

\begin{lemma}\label{L3.2}
 For every $q\in[1,\infty)$,
$\mathcal{L}^q(V)$ is a Banach space.
Moreover, it is reflexive for $q\in(1,\infty)$.
\end{lemma}

\begin{proof}
 Let us first consider the completeness of the space $\mathcal{L}^q $.  Let $\left\{\mathbf{u}_n\right\}$ be  a Cauchy sequence in $\mathcal{L}^q  $, where $\mathbf{u}_n=\left(u_{1,n}, u_{2,n}, \cdots,u_{i,n}, \cdots\right)$.  Then for any $\epsilon \in(0,1)$, there exists an integer $N>1$ such that for any $n, m > N$,
\begin{align}\label{3.2}
\left\|\mathbf{u}_n-\mathbf{u}_m\right\|_{\mathcal{L}^q }=\left(\sum_{x \in V} \mu(x)\left(\sum_{i=1}^{\infty}\left(u_{i, n}(x)-u_{i, m}(x)\right)^2\right)^{\frac{q}{2}}\right)^{\frac{1}{q}}<\epsilon.
\end{align}
For any fixed  $x_0 \in V$,  $\inf_{x \in V}\mu(x)>0$ and  \eqref{3.2}  imply that for any positive integer $M$, we get
\begin{align}\label{3.3}
\left(\sum_{i=1}^M\left(u_{i, n}(x_0)-u_{i, m}(x_0)\right)^2\right)^{\frac{1}{2}}< \mu_0^{-\frac{1}{q}}   \epsilon ,\quad \forall n,m>N.
\end{align}
Moreover,     $\left\{u_{i, n}(x_0)\right\}$ is a Cauchy sequence in $\mathbb{R}$ for all $i=1,2, \cdots$. Therefore, there exists  a real number $u_i(x_0)$  such that
  \begin{align}\label{u-n}
   \lim_{n\rightarrow\infty} u_{i, n}(x_0)=u_i(x_0), \quad \forall i=1,2, \cdots.
   \end{align}
By taking $n\rightarrow\infty$ in \eqref{3.3},
  we obtain
\begin{align*}
\left(\sum_{i=1}^M\left(u_{i}(x_0) -u_{i,m}(x_0)\right)^2\right)^{\frac{1}{2}}\leq \mu_0^{-\frac{1 }{ q}}   \epsilon,\quad \forall m>N.
\end{align*}
This together with \eqref{3.3}   leads to
\begin{align}\label{e5}
 \left(\sum_{i=1}^M\left(u_{i, n}(x_0)-u_i(x_0)\right)^2\right)^{\frac{1}{2}}
\nonumber\leq &\left(\sum_{i=1}^M\left(u_{i, n}(x_0)-u_{i, m}(x_0)\right)^2\right)^{\frac{1}{2}}+\left(\sum_{i=1}^M\left(u_{i, m}(x_0) -u_i(x_0) \right)^2\right)^{\frac{1}{2}} \\
 \leq & 2\mu_0^{-\frac{1}{q}}   \epsilon
\end{align}
for any $n > N$. Let $\mathbf{u}(x_0)=\left(u_1(x_0), u_2(x_0), \cdots,u_i(x_0), \cdots\right)$, where $u_i(x_0)$ is given by \eqref{u-n}.
Since    $\epsilon$ and $N$ are independent of $M$,
by taking $M\rightarrow \infty$ in \eqref{e5},
  we have
\begin{align}\label{3.4}
\left| \mathbf{u}_n(x_0)-\mathbf{u}(x_0)\right|
=\left(\sum_{i=1}^{\infty}\left(u_{i, n}(x_0)-u_i(x_0)\right)^2\right)^{\frac{1}{2}}
 \leq 2\mu_0^{-\frac{1}{q}}   \epsilon, \quad \forall n > N.
\end{align}
Since $|\mathbf{u}_n(x_0)|<\infty$, this also implies
that $|\mathbf{u}(x_0)|<\infty$.
 Noticing that \eqref{3.4} still holds if $x_0$ is replaced by any $x \in V$, we obtain
\begin{align*}
\lim_{n\rightarrow\infty}\left| \mathbf{u}_n(x)-\mathbf{u}(x)\right| =0,
\end{align*}
where $\mathbf u(x)$ is constructed in the same way as $\mathbf u(x_0)$ for each $x\in V$.
For $r>1$, write
$B_r(x_0)=\{y\in V:\mathrm d(x_0,y)<r\}$
for the ball centered at $x_0$ with radius $r$.
Since the graph is locally finite and connected, each $B_r(x_0)$ is finite, and these balls exhaust $V$ as $r\rightarrow\infty$.
Then we obtain
\begin{align*}
\lim_{n\rightarrow\infty}\left\|\mathbf{u}_n-\mathbf{u} \right\|_{\mathcal{L}^q(B_r(x_0))}
=\left(\sum_{x\in B_r(x_0)}  \mu(x) \lim_{n\rightarrow\infty}|\mathbf{u}_n(x)-\mathbf{u}(x)|^q \right)^{\frac{1}{q}}
=0;
\end{align*}
 namely   for the above $\epsilon$ and $N$, there exists a sufficiently large integer $n_r>N$     such that
\begin{align*}
\|\mathbf{u}_{n_r}-\mathbf{u} \|_{\mathcal{L}^q(B_r(x_0))}< \epsilon.
\end{align*}
This together with \eqref{3.2} leads to
\begin{align*}
\left\|\mathbf{u}_{n}-\mathbf{u}\right\|_{\mathcal{L}^q(B_r(x_0))}\leq
 \|\mathbf{u}_{n}-\mathbf{u}_{n_r} \|_{\mathcal{L}^q(B_r(x_0))}+
 \|\mathbf{u}_{n_r}-\mathbf{u} \|_{\mathcal{L}^q(B_r(x_0))}< 2\epsilon
\end{align*}
for any $n > N$.
Since $\epsilon$ and $N$ are independent of $r$, it holds that
\begin{align*}
\left\|\mathbf{u}_{n}-\mathbf{u}\right\|_{\mathcal{L}^q }\leq 2\epsilon, \quad \forall n > N.
\end{align*}
In particular, for any fixed $n>N$,
$$
\int_V|\mathbf{u}|^q\,d\mu
\leq
2^{q-1}\left(
\int_V|\mathbf{u}-\mathbf{u}_n|^q\,d\mu
+\int_V|\mathbf{u}_n|^q\,d\mu
\right)
<\infty,
$$
which leads to $\mathbf{u}\in\mathcal{L}^q$.
As a consequence, the  Cauchy sequence $\left\{\mathbf{u}_n\right\}$  converges to $\mathbf{u}$ in $\mathcal{L}^q $. This implies the completeness of $\mathcal{L}^q$, and thus $\mathcal{L}^q$ is a Banach space.

 It remains to prove reflexivity.
For $2\leq q<\infty$, the uniform convexity argument in
\cite[Proposition~3.2]{S-Y-Z-1} uses Hilbert space
inequalities that also hold for square-summable sequences.
It therefore applies to $\mathcal L^q(V)$, which is
uniformly convex and hence reflexive.
For $1<q<2$, set $q'=q/(q-1)>2$.
The same duality argument gives a linear isometric embedding
$$
T:\mathcal{L}^q
\rightarrow
\bigl(\mathcal{L}^{q'}\bigr)^*,
\quad
\langle T\mathbf{u},\mathbf{v}\rangle
=
\int_V\mathbf{u}\cdot\mathbf{v}\,d\mu.
$$
Since $\mathcal{L}^{q'}$ is reflexive, its dual is
reflexive.
Moreover, since $\mathcal{L}^q$ is complete and
$T$ is an isometry, $T(\mathcal{L}^q)$ is a closed
subspace of $\bigl(\mathcal{L}^{q'}\bigr)^*$.
Consequently, $\mathcal{L}^q$ is reflexive for $1<q<\infty$.
 \end{proof}

Using the spaces $L^p(V)$ and $\mathcal{L}^p(V)$,
we now establish Theorem~\ref{Sobolev1}(i).

\begin{lemma}\label{L3.3}
$W^{s,p}(V)$ is a reflexive Banach space.
\end{lemma}

\begin{proof}
By \cite{S-Y-Z-1} and Lemma~\ref{L3.2},
both $L^p(V)$ and $\mathcal{L}^p(V)$ are reflexive
Banach spaces for $p\in[2,\infty)$.
Let $\{u_n\}$ be a Cauchy sequence in $W^{s,p}(V)$.
Then $\{u_n\}$ and $\{\nabla^s u_n\}$ are Cauchy
sequences in $L^p(V)$ and $\mathcal{L}^p(V)$,
respectively.
Hence there exist $u\in L^p(V)$ and
$\mathbf{v}=(v_1,v_2,\cdots,v_i,\cdots)\in\mathcal{L}^p(V)$
such that
\begin{align}\label{e4}
u_n\rightarrow u \ \text{ in }\ L^p(V),
\quad
\nabla^s u_n\rightarrow\mathbf{v}
\ \text{ in }\ \mathcal{L}^p(V).
\end{align}
Since $\mu_0 >0$, we have
$$
|u(x)|
\leq \mu(x)^{-\frac{1}{p}}\|u\|_{L^p}
\leq \mu_0^{-\frac{1}{p}}\|u\|_{L^p},
\quad \forall x\in V.
$$
Thus $u\in L^\infty(V)$.
Moreover, for every $x\in V$,
$$
|u_n(x)-u(x)|
\leq \mu_0^{-\frac{1}{p}}\|u_n-u\|_{L^p}
\rightarrow0.
$$
This implies that  $u_n\rightarrow u$ pointwise  in  $V$ as $n \rightarrow\infty$.

 By H\"older's inequality and the convergence in
 \eqref{e4},
we obtain
   \begin{align}\label{3.7}
\int_V \nabla^s u_n(x) \cdot \boldsymbol{\varphi}(x) d \mu \rightarrow \int_V \mathbf{v}(x)\cdot \boldsymbol{ \varphi}(x) d \mu, \quad \forall \boldsymbol{\varphi} \in \mathcal{L}^{p'} ,
\end{align}
where   $p'=p/(p-1)$.
Fix $x_0\in V$ and $i\geq1$.
Let $\boldsymbol{\varphi}_i$ be the vector-valued
function satisfying   $\boldsymbol{\varphi}_i(x_0)=\mathbf{e}_i$ and    $\boldsymbol{\varphi}_i(x)=\mathbf{0}$ for $x\neq x_0$.
Since
$$\left\|\boldsymbol{\varphi}_i\right\|_{\mathcal{L}^{p'} }^{p'}
=\mu(x_0)<\infty,$$
 we have $\boldsymbol{\varphi}_i\in\mathcal{L}^{p'}(V)$. Taking $\boldsymbol{\varphi}=\boldsymbol{\varphi}_i$
in \eqref{3.7} and using the pointwise convergence
of $u_n$, we obtain
\begin{align*}
v_i(x_0)
=\sqrt{\frac{\widetilde{W}_s (x_0, x_i )}{2 \mu(x_0)}}\left(u(x_0)-u\left(x_i\right)\right).
\end{align*}
Since $i$ and $x_0$ are arbitrary, the definition of the fractional gradient in \eqref{e-7} gives
\begin{align*}
\mathbf{v}(x)=\nabla^su(x),  \quad \forall x\in V.
\end{align*}
From $u\in L^p(V)$ and
$\nabla^s u=\mathbf{v}\in\mathcal{L}^p(V)$,
we have $u\in W^{s,p}(V)$.
Moreover,
$$
\|u_n-u\|_{W^{s,p}}^p
=
\|u_n-u\|_{L^p}^p
+
\|\nabla^s u_n-\mathbf{v}\|_{\mathcal{L}^p}^p
\rightarrow0.
$$
Thus, $W^{s,p}(V)$ is complete and hence a Banach space.

Finally, equip
$X=L^p(V)\times\mathcal{L}^p(V)$
with the norm
$$ \|(u,\mathbf{v})\|_X=\left(\|u\|_{L^p}^p+\|\mathbf{v}\|_{\mathcal{L}^p}^p
\right)^{\frac{1}{p}}.$$
The map $T:W^{s,p}(V)\rightarrow X$ defined by $Tu=(u,\nabla^s u)$
is linear and isometric.
Since $W^{s,p}(V)$ is complete,
$T(W^{s,p}(V))$ is a closed subspace of $X$.
Both factors of $X$ are reflexive, so $X$ is reflexive.
Therefore, $T(W^{s,p}(V))$, and hence $W^{s,p}(V)$,
is reflexive.
\end{proof}

We next establish the strict inclusion and density asserted in Theorem~\ref{Sobolev1}(ii).

\begin{lemma}\label{L3.1}
 For any  $s\in (0,1)$ and $p\in [2,\infty)$,  $C_c(V)$ is a proper dense subspace of $W^{s,p}(V)$.
Consequently, $W^{s,p}_0(V)=W^{s,p}(V)$.
\end{lemma}

  \begin{proof}
  We first prove that $C_c(V)\subseteq W^{s,p}(V)$.  The conclusion is immediate if $u\equiv0$.
Let $u\in C_c(V)$ be nonzero.
Since $u$ has finite support, we have
$\int_V|u|^p\,d\mu<\infty$. It remains to verify that
\begin{align}\label{1.3}
\int_V|\nabla^s u|^p\,d\mu<\infty.
\end{align}
It follows from $u\in C_c(V)$ that the support $\mathrm{supp} (u)$ is finite.
We partition the set $V$  into two subsets, $\mathrm{supp} (u)$ and $V\backslash \mathrm{supp} (u)$.

 On the one hand,
for any fixed $x\in \mathrm{supp} (u)$,  we obtain
 \begin{align*}
 \sum_{y \in   \mathrm{supp} (u)\backslash \{x\}} W_s(x, y) (u(x)-u(y))^2  \leq 2\left(u^2(x)+ \max_{{y \in   \mathrm{supp} (u) }}u^2(y)\right)   \sum_{y \in  V\backslash \{x\}} W_s(x, y)
 \end{align*}
 and
  \begin{align*}
 \sum_{y \in  V\backslash \mathrm{supp} (u) } W_s(x, y)  (u(x)-u(y))^2
 \leq     u^2(x)   \sum_{y \in  V\backslash \{x\}} W_s(x, y).
\end{align*}
 Consequently, it follows from \eqref{kernel} that
 \begin{align}\label{e2}
  \nonumber&\int_{\mathrm{supp} (u)} |\nabla^su|^p   d\mu\\
\nonumber=&{2^{-\frac{p}{2}}}\sum_{x\in \mathrm{supp} (u)} \left(\sum_{y \in V \backslash \{x\}} W_s(x, y) (u(x)-u(y))^2\right)^{\frac{p}{2}}\mu^{1-\frac{p}{2}} (x)\\
\nonumber\leq &{2^{-\frac{p}{2}}}
\max_{x \in   \mathrm{supp} (u) }\mu^{1-\frac{p}{2}} (x)
\sum_{x\in \mathrm{supp} (u)} \left(3u^2(x)+ 2\max_{{y \in   \mathrm{supp} (u) }}u^2(y)\right) ^{\frac{p}{2}}\left(   \sum_{y \in  V\backslash \{x\}} W_s(x, y) \right)^{\frac{p}{2}} \\
 <& \infty.
\end{align}
On the other hand, by Tonelli's theorem, \eqref{symmetry} and \eqref{kernel}, we have
\begin{align}\label{W-1}
\sum_{x\in V\setminus\mathrm{supp}(u)}
\sum_{y\in\mathrm{supp}(u)}W_s(x,y)
&=
\sum_{y\in\mathrm{supp}(u)}
\sum_{x\in V\setminus\mathrm{supp}(u)}W_s(y,x)
\nonumber\\
&\leq
\sum_{y\in\mathrm{supp}(u)}
\sum_{x\in V\setminus\{y\}}W_s(y,x)
<\infty.
\end{align}
 Since $p\in [2,\infty)$  and $\mu_0  >0$, we have  $0< \mu(x)^{1- {p}/{2}} \leq  \mu_0^{1- {p}/{2}}$ for all $x\in V$,   and
$\sum_n a_n^{p/2}\leq(\sum_n a_n)^{p/2}$
for every nonnegative sequence $\{a_n\}$.
 Then we obtain
    \begin{align*}
 \int_{V\backslash \mathrm{supp} (u)} |\nabla^su|^p   d\mu
 =&   2^{-\frac{p}{2}} \sum_{x\in V\backslash \mathrm{supp} (u)}  \mu(x) ^{1-\frac{p}{2}}
 \left(\sum_{y \in   \mathrm{supp} (u)} W_s(x, y) u^2(y) \right)^{\frac{p}{2}}   \\
\leq & 2^{-\frac{p}{2}}  \mu_0 ^{1-\frac{p}{2}}  \max_{{y \in   \mathrm{supp} (u) }}|u|^p(y)  \sum_{x\in V\backslash \mathrm{supp} (u)} \left( \sum_{y \in   \mathrm{supp} (u)} W_s(x, y)   \right)^{\frac{p}{2}}\\
\leq& 2^{-\frac{p}{2}}  \mu_0 ^{1-\frac{p}{2}}  \max_{{y \in   \mathrm{supp} (u) }}|u|^p(y)  \left(\sum_{x\in V\backslash \mathrm{supp} (u)}  \sum_{y \in   \mathrm{supp} (u)} W_s(x, y)   \right)^{\frac{p}{2}}.
\end{align*}
According to this and \eqref{W-1}, we get
   \begin{align}\label{e3}
  \int_{V\backslash \mathrm{supp} (u)} |\nabla^su|^p   d\mu<\infty.
  \end{align}
  Therefore, \eqref{1.3} follows from  \eqref{e2} and \eqref{e3},  which leads to $ C_c(V)   \subseteq  W^{s,p}(V)$.

 We next show that the inclusion is strict.
Write $V=\{x_1,x_2,\cdots,x_i,\cdots\}$, and let
$\varphi_j$ be the indicator function of $x_j$.
By the inclusion proved above,
$\varphi_j\in W^{s,p}(V)$.
Set
\[
a_j=\frac{2^{-j}}{1+\|\varphi_j\|_{W^{s,p} }}\ \text{ and }
\
u_N=\sum_{j=1}^N a_j\varphi_j.
\]
For $N>M$, we have
\[
\|u_N-u_M\|_{W^{s,p} }
\leq
\sum_{j=M+1}^N a_j\|\varphi_j\|_{W^{s,p} }
\leq
\sum_{j=M+1}^N2^{-j}.
\]
Thus $\{u_N\}$ is a Cauchy sequence in $W^{s,p}(V)$.
By Lemma~\ref{L3.3}, there exists $u\in W^{s,p}(V)$
such that $u_N\rightarrow u$ in $W^{s,p}(V)$.
For each fixed $j$, the estimate
\[
|u_N(x_j)-u(x_j)|
\leq
\mu(x_j)^{-\frac{1}{p}}\|u_N-u\|_{L^p }
\leq
\mu(x_j)^{-\frac{1}{p}}\|u_N-u\|_{W^{s,p} }
\]
implies that $u(x_j)=a_j>0$.
Hence $\operatorname{supp}(u)=V$ is infinite, so
$u\notin C_c(V)$.
This proves that $C_c(V)\subsetneq W^{s,p}(V)$.

 Finally, we prove that $C_c(V)$ is dense in $W^{s,p}(V)$. Fix $u\in W^{s,p}(V)$ and define
$$
u_n=\operatorname{sgn}(u)\left(|u|-\frac1n\right)^+.
$$
For every finite subset
$A\subseteq \mathrm{supp}(u_n)= \{x\in V:|u(x)|>1/n\}$, we have
\[
\mu_0n^{-p}\#A
\leq
\sum_{x\in A}|u(x)|^p\mu(x)
\leq
\|u\|_{L^p }^p.
\]
Since $ \inf_{x\in V}\mu(x)>0$ and $u\in L^p(V)$,
the set $ \mathrm{supp}(u_n)$ is finite, and
thus $u_n\in C_c(V)$.

 Set
$$
\varphi_n =u-u_n
=\max\left\{-\frac1n,\min\left\{u,\frac1n\right\}\right\}.
$$
Then $\varphi_n\rightarrow0$ pointwise, $|\varphi_n|\leq|u|$, and
$$
|\varphi_n(x)-\varphi_n(y)|\leq|u(x)-u(y)|,
\quad \forall x,y\in V.
$$
For each fixed $x\in V$, the dominated convergence
theorem applied to the kernel sum defining the
fractional gradient gives
$$
|\nabla^s \varphi_n(x)|\rightarrow0,
\quad
|\nabla^s \varphi_n(x)|\leq|\nabla^s u(x)|.
$$
Applying the dominated convergence theorem once more,
we obtain
$$
\|u-u_n\|_{W^{s,p}}^p
=
\int_V\left(|\varphi_n|^p+|\nabla^s \varphi_n|^p\right)d\mu
\rightarrow0.
$$
Hence $C_c(V)$ is dense in $W^{s,p}(V)$.
Since $W^{s,p}(V)$ is complete by Lemma~\ref{L3.3},
its identification with the completion of $C_c(V)$
yields $W^{s,p}_0(V)=W^{s,p}(V)$.
\end{proof}

We next prove the continuous embeddings stated
in Theorem~\ref{Sobolev1}(iii).

\begin{lemma}\label{L3.4}
$W^{s,p}(V)$ is continuously embedded in $L^q(V)$
for every $q\in[p,\infty]$.
\end{lemma}

\begin{proof}
For any $u\in W^{s,p}(V)$, we have
$$
\|u\|_{L^\infty}
\leq \mu_0^{-\frac{1}{p}}\|u\|_{L^p}
\leq \mu_0^{-\frac{1}{p}}\|u\|_{W^{s,p}}.
$$
For $p<q<\infty$, it follows that
$$
\|u\|_{L^q}
 \leq
\|u\|_{L^\infty }^{1-\frac{p}{q} }
\|u\|_{L^p}^{\frac{p}{q}} \leq
\mu_0^{\frac{1}{q}-\frac{1}{p}}\|u\|_{W^{s,p}}.
$$
Together with
$\|u\|_{L^p}\leq\|u\|_{W^{s,p}}$,
these estimates prove this lemma.
\end{proof}

 Furthermore, we verify the properties of the positive and negative parts of $ u $, namely Theorem \ref{Sobolev1}(iv).
  For any $u\in  L^\infty(V)$,  it is clear that $u^+\geq 0$, $u^-\leq 0$, $u=u^++u^-$  and  $u^+u^-=0$.
Fixing $x\in V$ and using the definition in \eqref{gra}, we compute
\begin{align}\label{13}
\nonumber  {\nabla ^su^+\nabla ^su^-}(x)   &  =   \frac{1}{2\mu(x)} \sum_{y \in V \backslash \{x\}} W_s(x, y) \left(u^+(x)-u^+(y)\right)\left(u^-(x)-u^-(y)\right)\\
\nonumber &  =-   \frac{ 1}{2\mu(x)} \sum_{y \in V \backslash \{x\}} W_s(x, y) \left(  u^+(x)u^-(y)+  u^+(y)u^-(x)\right)\\
&\geq0.
\end{align}
 This inequality will be helpful in proving the following lemmas.

\begin{lemma}\label{L3.5}
If $u\in W^{s,p}(V)$, then $u^+$, $u^-$ and  $|u|$   also belong to $W^{s,p}(V)$.
\end{lemma}

\begin{proof}
For any $u\in W^{s,p}(V)$, it follows from  $u^+=u- u^-$,   \eqref{13}  and the  Minkowski  inequality that
\begin{align*}
 \|u^+\|_{W^{s,p}}  
  \leq &\left(\int_{ V} |\nabla^s u^+|^p    d\mu\right)^{\frac{1}{p}}+\left(\int_{ V} (u^+)^p  d\mu\right)^{\frac{1}{p}}\\
    =&\left(\int_{ V} \left( |\nabla^s (u- u^-) |^2\right)^ {\frac{p}{2} }d\mu\right)^{\frac{1}{p}}+\left(\int_{ V} \left(\left(u- u^-\right)^2\right)^ {\frac{p}{2} } d\mu\right)^{\frac{1}{p}}\\
=&\left(\int_{ V} \left( |\nabla^s  u|^2-| \nabla^su^-  |^2-2 \nabla^s u^+\nabla^su^- \right)^ {\frac{p}{2} }d\mu\right)^{\frac{1}{p}}+\left(\int_{ V} \left( u^2-( u^-)^2 \right)^ {\frac{p}{2} } d\mu\right)^{\frac{1}{p}}\\
\leq&\left(\int_{ V} |\nabla^s  u|^pd\mu\right)^{\frac{1}{p}}+\left(\int_{ V}  | u|^p   d\mu\right)^{\frac{1}{p}}\\
\leq & 2\|u\|_{W^{s,p}} .
\end{align*}
This  leads to $u^+ \in W^{s,p}(V)$.  A similar argument establishes  $u^- \in W^{s,p}(V)$.
Since $W^{s,p}(V)$ is a linear space and $|u|=u^+-u^-$, we have $|u|\in W^{s,p}(V)$, and  this completes the proof.
\end{proof}

\begin{proof}[\textbf{\textit{Proof of Theorem  \ref{Sobolev1}}}]
Statements (i)--(iv) follow from
Lemmas~\ref{L3.3}, \ref{L3.1}, \ref{L3.4}
and~\ref{L3.5}, respectively.
\end{proof}

\section{Fractional  $p$-Laplace equation}\label{S4}

\subsection{Preliminaries}

We first establish the following properties of
$\mathcal{H}_{s,p}$, defined by \eqref{H}.
\begin{lemma}\label{Le4-1}
 The space $\mathcal{H}_{s,p}$ is a reflexive Banach space,
and is continuously and compactly embedded in $L^q(V)$
for every $q\in[p,\infty]$.
More precisely, if $\{u_n\}$ is bounded in $\mathcal{H}_{s,p}$, then   there exists $u\in\mathcal{H}_{s,p}$ such that, up to a subsequence, as $n\rightarrow\infty$,
 \begin{align*}\left\{
        \begin{aligned}
			&u_n\rightharpoonup u \ \text{ weakly in }  \mathcal{H}_{s,p},\\
			&u_n\rightarrow u \ \text{ strongly in } L^q(V) \   \text{ for any }  q\in [p,\infty],\\
            &u_n\rightarrow u\ \text{  pointwise  in  }   V.
        \end{aligned}\right.
      	\end{align*}
\end{lemma}

\begin{proof}
By definition, ${\mathcal{H}_{s,p}}$ is a Banach space.
To prove reflexivity, let
$X=L^p(V)\times\mathcal{L}^p(V)$
be equipped with the norm
$$ \|(u,\mathbf{v})\|_X=\left(\|u\|_{L^p}^p+\|\mathbf{v}\|_{\mathcal{L}^p}^p
\right)^{\frac{1}{p}}.$$
Then $X$ is a reflexive Banach space by Lemma \ref{L3.2} and the
reflexivity of $L^p(V)$.
For any $u\in C_c(V)$, define
$$
T_0u=\left(h^{\frac{1}{p}}u,\nabla^s u\right).
$$
This map is linear and preserves the norm, because
$$
\|T_0u\|_X^p
=
\int_V\bigl(h|u|^p+|\nabla^s u|^p\bigr)\,d\mu
=
\|u\|_{\mathcal{H}_{s,p}}^p.
$$
 Since $\mathcal{H}_{s,p}$ is the completion of
$C_c(V)$ under the norm $\|\cdot\|_{\mathcal{H}_{s,p}}$,
$T_0$ extends uniquely to a linear isometry
$T:\mathcal{H}_{s,p}\rightarrow X $
satisfying
$$
\|Tu\|_X=\|u\|_{\mathcal{H}_{s,p}},
\quad \forall u\in\mathcal{H}_{s,p}.
$$
Since $\mathcal H_{s,p}$ is complete and $T$ is an isometry,
$T(\mathcal H_{s,p})$ is a closed subspace of $X$.
As $X$ is reflexive, so are $T(\mathcal H_{s,p})$
and $\mathcal H_{s,p}$.

 Moreover, \eqref{H-1} implies that, for every $u\in C_c(V)$,
\begin{align*}
\|u\|_{W^{s,p}}^p
&=\int_V\bigl(|\nabla^s u|^p+|u|^p\bigr)\,d\mu\\
&\leq\int_V|\nabla^s u|^p\,d\mu
 +\frac{1}{h_0}\int_Vh|u|^p\,d\mu\\
&\leq\left(1+\frac{1}{h_0}\right)
       \|u\|_{\mathcal{H}_{s,p}}^p.
\end{align*}
By completion, this estimate defines a continuous linear map
from $\mathcal H_{s,p}$ to $W^{s,p}(V)$.
If $\{u_n\}\subseteq C_c(V)$ is Cauchy in $\mathcal H_{s,p}$
and converges to $u$ in $W^{s,p}(V)$, then
$h^{1/p}u_n$ converges in $L^p(V)$ to $h^{1/p}u$,
as follows from convergence at each vertex.
Together with the convergence of $\nabla^s u_n$ to $\nabla^s u$,
this shows that the map is injective.
Consequently, $\mathcal{H}_{s,p}$ is continuously
embedded in $W^{s,p}(V)$.
Theorem \ref{Sobolev1}(iii) then yields its continuous
embedding into $L^q(V)$ for every $q\in[p,\infty]$.

 To prove compactness, let $\{u_n\}$ be a bounded sequence in
${\mathcal{H}_{s,p}}$.
Since $\mathcal H_{s,p}$ is reflexive,
there exist a subsequence,
still denoted by $\{u_n\}$, and $u\in\mathcal H_{s,p}$ such that
 $u_n\rightharpoonup u$ weakly in ${\mathcal{H}_{s,p}}$ as ${n\rightarrow\infty}$.
 Furthermore,
for any fixed vertex $x\in V$ and every $v\in \mathcal{H}_{s,p} $, the estimate
$$|v(x)|^p\leq \frac{1}{\mu(x)h(x)}     \int_V h|v|^p\,d\mu\leq \frac{1}{\mu(x)h(x)}    \|v\|_{ \mathcal{H}_{s,p} }^p $$
shows that point evaluation is a continuous linear functional
on $\mathcal{H}_{s,p}$.
Hence $u_n$ converges pointwise to $u$ on $V$ as $n\to\infty$.

 We next prove that $u_n\rightarrow u$ strongly in $L^p(V)$.
By the boundedness of $\{u_n\}$, there exists a constant $M>0$ 
such that
	\begin{align}\label{9}
\left(\int_{V}h|u_n-u|^p  d\mu\right)^{\frac{1}{p}}
\leq \|u_n -u \|_{\mathcal{H}_{s,p}}
\leq \|u_n \|_{\mathcal{H}_{s,p}} +\| u \|_{\mathcal{H}_{s,p}}  \leq M .
	\end{align}
Fix a vertex $x_0\in V$.
For every $\epsilon>0$, assumption \eqref{H-2} yields
a sufficiently large $r\geq 1$ 
 such that
\begin{align}\label{10-1}
h(x)\geq \frac{M^p}{\epsilon},
\quad \forall x\in V\setminus B_r(x_0).
\end{align}
It follows from \eqref{9} and \eqref{10-1} that
	\begin{align}\label{11-1}
	\int_{V\backslash B_r(x_0)} |u_n-u|^pd\mu
 \leq   \frac{\epsilon}{M^p}
 \int_{V} h  |u_n-u|^pd\mu
 \leq   \epsilon.
	\end{align}	
 Since the graph is locally finite, $B_r(x_0)$ contains only
finitely many vertices. Hence the pointwise convergence implies
$$
\lim_{n\rightarrow\infty}
\int_{B_r(x_0)}|u_n-u|^p\,d\mu=0.
$$
Together with \eqref{11-1}, this gives
$$
\limsup_{n\rightarrow\infty}
\int_V|u_n-u|^p\,d\mu\leq\epsilon.
$$
Since $\epsilon>0$ is arbitrary, we conclude that
$u_n\rightarrow u$ strongly in $L^p(V)$.

 Moreover, we derive that  $u_n$ converges strongly to $ u$  in ${L^\infty}(V)$  by noting that
\begin{align*}
\|u_n-u\|_{L^\infty}^p\leq\frac{1}{\mu_0}\int_{V}|u_n-u|^pd\mu\rightarrow0.
\end{align*}
 For every $p<q<\infty$, we also have
$$
\|u_n-u\|_{L^q}^q
\leq
\|u_n-u\|_{L^\infty}^{q-p}
\|u_n-u\|_{L^p}^p
\rightarrow 0.
$$
Then $u_n\rightarrow u$ strongly in $L^q(V)$ for every
$q\in[p,\infty]$. This completes the proof.
\end{proof}

Secondly, we establish a positivity result for weak solutions.

\begin{lemma}\label{L4.2}
Suppose $ f(x, y) \equiv 0$   for any $ y\leq 0$.
If  $u\in\mathcal H_{s,p}(V)$ is a nontrivial weak solution of \eqref{eS1}, then it is also
 a strictly positive solution  on $V$.
\end{lemma}

\begin{proof}
Let $u\in\mathcal H_{s,p}(V)$ be a nontrivial weak
solution of \eqref{eS1}.
 We first verify that $u^-\in\mathcal H_{s,p}(V)$.
Choose $\{u_n\}\subseteq C_c(V)$ such that
$u_n\rightarrow u$ in $\mathcal H_{s,p}(V)$.
Then $u_n^-\in C_c(V)$, and the inequalities
$|u_n^-|\leq|u_n|$ and
$|\nabla^s u_n^-|\leq|\nabla^s u_n|$ give
\begin{align*}
\|u_n^-\|_{\mathcal H_{s,p}}^p
=
\int_V\left(
|\nabla^s u_n^-|^p+h|u_n^-|^p
\right)\,d\mu
\leq
\|u_n\|_{\mathcal H_{s,p}}^p.
\end{align*}
Thus $\{u_n^-\}$ is bounded in $\mathcal H_{s,p}(V)$.
By reflexivity,
there exists $v\in\mathcal H_{s,p}(V)$ such that
$u_n^-\rightharpoonup v$ in $\mathcal H_{s,p}(V)$,
and then  $u_n^-(x)\rightarrow v(x)$ for every $x\in V$.
On the other hand, $u_n(x)\rightarrow u(x)$ implies
$u_n^-(x)\rightarrow u^-(x)$.
Hence $v=u^-$, and therefore
$u^-\in\mathcal H_{s,p}(V)$.

 Taking $\varphi=u^-$ in \eqref{weak-s} and using
\eqref{H-1} and \eqref{13}, we obtain
$u^-\equiv0$ on $V$. Thus $u\geq0$ on $V$.
Suppose that $u(x_0)=0$ for some $x_0\in V$.
Since $u$ is a pointwise solution of \eqref{eS1}
and $f(x_0,0)=0$, we have
$(-\Delta)^s_pu(x_0)=0$, that is,
\begin{align}\label{1}
\frac{1}{2\mu(x_0)}
\sum_{y\in V\setminus\{x_0\}}W_s(x_0,y)
\left(
|\nabla^s u|^{p-2}(x_0)
+|\nabla^s u|^{p-2}(y)
\right)u(y)=0.
\end{align}
Since $u\geq0$ is nontrivial,
\eqref{symmetry} and \eqref{N} imply
that $|\nabla^s u|(x_0)>0$.
Consequently, every coefficient of $u(y)$ in \eqref{1} is strictly positive.
As $u(y)\geq0$, it follows that $u(y)=0$ for
every $y\neq x_0$.
Together with $u(x_0)=0$, this gives $u\equiv0$ on $V$,
contradicting the nontriviality of $u$.
Therefore, $u(x)>0$ for every $x\in V$.
 \end{proof}

We next establish a vector monotonicity inequality.

\begin{lemma}\label{L4.5}
Let  $p\in [2,\infty)$.
 For any  vectors $ \mathbf{a}$, $ \mathbf{b}\in \mathbb{R}^{\infty}$ satisfying  $ |\mathbf{a}|<\infty$ and $ |\mathbf{b}|<\infty$,  we have
\begin{align}\label{inq-2}
   \left ( | \mathbf{a}|^{p-2}   \mathbf{a} -|\mathbf{b}|^{p-2}  \mathbf{b}  \right) \cdot(\mathbf{a}-  \mathbf{b})
    \geq
 2^{2-p} | \mathbf{a}-\mathbf{b} |^p .
\end{align}
 Moreover,   the constant $2^{2-p}$ is optimal.
\end{lemma}

\begin{proof}
 To prove the inequality \eqref{inq-2}, we adapt the finite-dimensional argument
in \cite[Lemma~A.0.5]{Peral} to the present infinite-dimensional setting.
We include the details for completeness.
 For $p=2$, the inequality is an identity. Hence, we assume that $p>2$.
Define
\begin{align*}
 \phi(t, m)=\frac{1+t^p- (t^{p-1}+t ) m}{\left (1+t^2-2 t m \right)^{\frac{p}{2}}}
\end{align*}
for any $t\in (0,1)$ and $m\in[-1,1]$.

 We claim that
 \begin{align}\label{e1}
 \phi(t, m) \geq   2^{2-p}.
 \end{align}
Fix $t\in(0,1)$ and $m\in[-1,1]$, and set $
D=1+t^2-2tm$.
It is clear that
$$
0<(1-t)^2\leq D\leq(1+t)^2.
$$
A direct calculation gives
\begin{align*}
1+t^p-(t^{p-1}+t)m
&=
\frac{1+t^{p-1}}{1+t}\,D
+\frac{t(1-t)(1-t^{p-2})(1+m)}{1+t}\\
&\geq\frac{1+t^{p-1}}{1+t}\,D.
\end{align*}
Dividing by $D^{p/2}$ and using
$D\leq(1+t)^2$ and $p>2$, we obtain
$$
\phi(t,m)
\geq
\frac{1+t^{p-1}}{(1+t)D^{\frac{p}{2}-1} }
\geq
\frac{1+t^{p-1}}{(1+t)^{p-1}}.
$$
Finally, the convexity of $r\mapsto r^{p-1}$ on
$[0,\infty)$ yields
$$
1+t^{p-1}
\geq
2^{2-p}(1+t)^{p-1}.
$$
Combining the preceding two inequalities, we prove  \eqref{e1}.

If $|\mathbf a|=0$ or $|\mathbf b|=0$, then
\eqref{inq-2} is immediate.
If $|\mathbf a|=|\mathbf b|=r>0$, then
$|\mathbf a-\mathbf b|\leq2r$, and hence
\begin{align*}
 \bigl(|\mathbf a|^{p-2}\mathbf a
      -|\mathbf b|^{p-2}\mathbf b\bigr)
 \cdot(\mathbf a-\mathbf b)
 =r^{p-2}|\mathbf a-\mathbf b|^2
\geq2^{2-p}|\mathbf a-\mathbf b|^p.
\end{align*}
If $|\mathbf{a}|   >  |\mathbf{b}|   >0$, set
$$
t=\frac{|\mathbf b|}{|\mathbf a|},
\quad
m=\frac{\mathbf a\cdot\mathbf b}
        {|\mathbf a|\,|\mathbf b|}.
$$
Then $t\in(0,1)$ and $|m|\leq1$ by the
Cauchy--Schwarz inequality. Consequently,
$$
\phi(t,m)=
\frac{
\bigl(|\mathbf a|^{p-2}\mathbf a
      -|\mathbf b|^{p-2}\mathbf b\bigr)
\cdot(\mathbf a-\mathbf b)}
{|\mathbf a-\mathbf b|^p}
\geq 2^{2-p}.
$$
The case $|\mathbf b|>|\mathbf a|>0$ follows by interchanging
$\mathbf a$ and $\mathbf b$.
Thus, \eqref{inq-2} holds in all cases.

 Finally, taking $\mathbf b=-\mathbf a\neq\mathbf0$, we have
$$
\left(|\mathbf a|^{p-2}\mathbf a
      -|\mathbf b|^{p-2}\mathbf b\right)
\cdot(\mathbf a-\mathbf b)
=4|\mathbf a|^p
=2^{2-p}|\mathbf a-\mathbf b|^p.
$$
Thus, the constant $2^{2-p}$ is optimal.
\end{proof}

\subsection{Mountain-pass type solutions}\label{sub2}

We first establish the mountain-pass geometry of $E_{s,p}$.

\begin{lemma}\label{L4.3}
There exist a  nonnegative function
$u\in\mathcal H_{s,p}(V)$ and a constant $r>0$ such that
\begin{align*}
\lim_{t\rightarrow+\infty}E_{s,p}(tu)=-\infty,
\quad
\inf_{\|v\|_{\mathcal H_{s,p}}=r}E_{s,p}(v)>0.
\end{align*}
\end{lemma}

\begin{proof}
 Fix $x_0\in V$ and define $u$ by $u(x_0)=1$ and
$u(x)=0$ for $x\neq x_0$.
By Lemma~\ref{L3.1} and
\begin{align*}
\int_Vh|u|^p\,d\mu=h(x_0)\mu(x_0)<\infty,
\end{align*}
we have $u\in C_c(V)\subseteq\mathcal H_{s,p}(V)$.

 By $(\mathrm{A}_3)$, we obtain   $ F(x, y)>0$ and
\begin{align*}
 \frac{\partial}{\partial y} \log \frac{F(x, y)} {y^{\alpha}}=\frac{yf(x,y)-\alpha F(x,y)}{y F(x,y)}   \geq 0, \quad \forall x\in V, \ y>0.
\end{align*}
Thus, for every fixed $x\in V$, the function $y\mapsto F(x,y)/y^\alpha$ is nondecreasing on $(0,\infty)$, and hence
\begin{align}\label{eq:F-scaling}
F(x,ty)\ge t^\alpha F(x,y),
\quad \forall x\in V,\, y\ge0,\, t\ge1.
\end{align}
 Consequently, for $t\geq1$,
\begin{align*}
E_{s,p}(tu)
&=
\frac{t^p}{p}\|u\|_{\mathcal H_{s,p}}^p
-\mu(x_0)F(x_0,t)\\
&\leq
\frac{t^p}{p}\|u\|_{\mathcal H_{s,p}}^p
-\mu(x_0)F(x_0,1)t^\alpha.
\end{align*}
Since $F(x_0,1)>0$ and $\alpha>p$, we obtain
$E_{s,p}(tu)\rightarrow-\infty$ as $t\rightarrow+\infty$.

 Next, by the uniformity in $(\mathrm{A}_4)$ and
$f(x,0)=0$, there exist $\epsilon\in(0,\lambda_p)$
and $\delta>0$ such that
\begin{align*}
f(x,y)\leq(\lambda_p-\epsilon)y^{p-1},
\quad \forall x\in V,\, y\in [0,\delta].
\end{align*}
Integrating with respect to $y$, we obtain
\begin{align}\label{F-small-bound}
F(x,y)\leq\frac{\lambda_p-\epsilon}{p}y^p,
\quad \forall x\in V,\, y\in [0,\delta].
\end{align}
By Lemma~\ref{Le4-1}, \eqref{mu} and \eqref{H-1}, we have
\begin{align*}
\|v\|_{L^\infty}
\leq (\mu_0h_0)^{-\frac{1}{p}}\|v\|_{\mathcal H_{s,p}},
\quad \forall v\in\mathcal H_{s,p}(V).
\end{align*}
Set $r= (\mu_0h_0)^{1/p}\delta/2$.
For every $v\in\mathcal H_{s,p}(V)$ with
$\|v\|_{\mathcal H_{s,p}}=r$, we have
$0\leq v^+(x)\leq\delta/2$ for all $x\in V$.
Thus, by \eqref{F-small-bound} and the definition of
$\lambda_p$,
\begin{align*}
E_{s,p}(v)
&\geq
\frac1p\|v\|_{\mathcal H_{s,p}}^p
-\frac{\lambda_p-\epsilon}{p}\|v^+\|_{L^p}^p\\
&\geq
\frac1p\|v\|_{\mathcal H_{s,p}}^p
-\frac{\lambda_p-\epsilon}{p\lambda_p}
\|v\|_{\mathcal H_{s,p}}^p\\
&=
\frac{\epsilon}{p\lambda_p}r^p.
\end{align*}
Therefore, we obtain
\begin{align*}
\inf_{\|v\|_{\mathcal H_{s,p}}=r}E_{s,p}(v)
\geq\frac{\epsilon}{p\lambda_p}r^p>0,
\end{align*}
which completes the proof.
\end{proof}

We next prove that  $E_{s,p} $ satisfies the  Palais--Smale condition.

\begin{lemma} \label{L3.4-2}
For any $c \in \mathbb{R}$,  the functional $E_{s,p} $ satisfies the $(\mathrm{PS})_c$   condition. Namely, if $\{u_n\} \subseteq \mathcal{H}_{s,p}$ satisfies $E_{s,p} \left(u_n\right) \rightarrow c$ and
 $ \|E_{s,p}'(u_n)\|_{(\mathcal H_{s,p}(V))^*}\rightarrow 0$ as $n\rightarrow\infty$,
then there exist a subsequence, still denoted by
$\{u_n\}$, and $u\in\mathcal H_{s,p}(V)$ such that
$u_n\to u$ strongly in $\mathcal H_{s,p}(V)$.
\end{lemma}

\begin{proof}
Since $f(x,0)=0$, we have
$f(x,u_n^+)u_n=f(x,u_n^+)u_n^+$.
Thus $(\mathrm{A}_3)$ gives
\begin{align*}
\left(\frac1p-\frac1\alpha\right)
\|u_n\|_{\mathcal H_{s,p}}^p
&\leq
E_{s,p}(u_n)
-\frac1\alpha\langle E'_{s,p}(u_n),u_n\rangle\\
&\leq
|c|+1+
\frac1\alpha\|E'_{s,p}(u_n)\|_{(\mathcal H_{s,p})^*}
\|u_n\|_{\mathcal H_{s,p}}
\end{align*}
for all sufficiently large $n$.
Since $\alpha>p>1$ and
$\|E'_{s,p}(u_n)\|_{(\mathcal H_{s,p})^*}\rightarrow0$,
the sequence $\{u_n\}$ is bounded in $\mathcal H_{s,p}$.

 By Lemma~\ref{Le4-1}, after passing to a subsequence,
still denoted by $\{u_n\}$,  there exists
$u\in\mathcal{H}_{s,p}$ such that
$u_n\rightharpoonup u$ weakly in $\mathcal{H}_{s,p}$ and
$u_n\rightarrow u$ strongly in $L^q$ for every $q\in[p,\infty]$.
In particular, $\{u_n\}$ is bounded in $L^\infty$
by a constant $M>0$.
 Assumptions $(\mathrm{A}_1)$--$(\mathrm{A}_4)$ give
a constant $K_M>0$ such that
\begin{align}\label{5}
|f(x,y)|\leq K_My^{p-1},
\quad \forall x\in V, \,  y\in  [0, M].
\end{align}
By \eqref{5} and H\"older's inequality,
\begin{align}\label{e21}
\left|\int_Vf(x,u_n^+)(u_n-u)\,d\mu\right|
&\leq
K_M\|u_n\|_{L^p}^{p-1}\|u_n-u\|_{L^p}
=o_n(1),
\end{align}
since $\{u_n\}$ is bounded in $L^p(V)$ and
$u_n\to u$ strongly in $L^p(V)$.
Furthermore, since $\{u_n-u\}$ is bounded in
$\mathcal H_{s,p}(V)$ and
$\|E_{s,p}'(u_n)\|_{(\mathcal H_{s,p}(V))^*}\to0$,
we have
$$
\left|\left\langle E'_{s,p}(u_n),u_n-u\right\rangle\right|
\leq
\|E'_{s,p}(u_n)\|_{(\mathcal{H}_{s,p})^*}
\|u_n-u\|_{\mathcal{H}_{s,p}}
=o_n(1).
$$
Combining this with
\eqref{E-C1-derivative} and \eqref{e21}, we obtain
\begin{align}\label{3}
\int_V \Big(
|\nabla^s u_n|^{p-2}\nabla^s u_n
\cdot\nabla^s(u_n-u)
+h|u_n|^{p-2}u_n(u_n-u)
\Big)\,d\mu
=o_n(1).
\end{align}
Define
$$
L_u(\varphi)
 =\int_V\Big(
|\nabla^s u|^{p-2}\nabla^s u\cdot\nabla^s\varphi
+h|u|^{p-2}u\varphi
\Big)\,d\mu,
\quad \forall \varphi\in\mathcal{H}_{s,p}.
$$
By H\"older's inequality,
$$
|L_u(\varphi)|
\leq
\|u\|_{\mathcal{H}_{s,p}}^{p-1}
\|\varphi\|_{\mathcal{H}_{s,p}}.
$$
Thus, $L_u$ is a continuous linear functional on
$\mathcal{H}_{s,p}$.
Since $u_n\rightharpoonup u$ weakly in $\mathcal{H}_{s,p}$,
we have
\begin{align}\label{4}
\int_V\Big(
|\nabla^s u|^{p-2}\nabla^s u
\cdot\nabla^s(u_n-u)
+h|u|^{p-2}u(u_n-u)
\Big)\,d\mu
=o_n(1).
\end{align}

 Finally, applying Lemma~\ref{L4.5} pointwise to the fractional gradients and to the scalar terms,
with real numbers identified with vectors of the form
$(a,0,\cdots,0,\cdots)$, and using $h>0$, we obtain
\begin{align*}
\|u_n-u\|_{\mathcal{H}_{s,p}}^p
&=\int_V\left(
|\nabla^s(u_n-u)|^p+h|u_n-u|^p
\right)d\mu\\
&\leq 2^{p-2}
\int_V
\left(
|\nabla^s u_n|^{p-2}\nabla^s u_n
-|\nabla^s u|^{p-2}\nabla^s u
\right)
\nabla^s(u_n-u)\,d\mu\\
&\quad+2^{p-2}
\int_V h
\left(
|u_n|^{p-2}u_n-|u|^{p-2}u
\right)(u_n-u)\,d\mu\\
&=o_n(1),
\end{align*}
where the last equality follows from \eqref{3} and \eqref{4}.
Consequently, $u_n\to u$ strongly in $\mathcal H_{s,p}(V)$.
This proves the $(\mathrm{PS})_c$ condition.
\end{proof}

\begin{proof}[\textbf{Proof of Theorem  \ref{T3}}]
Lemmas~\ref{L4.3} and~\ref{L3.4-2} establish
the mountain-pass geometry and the Palais--Smale
condition, respectively.
We apply the mountain-pass theorem in \cite{A-R}.
First, Lemma  \ref{L4.3}    ensures that there exists a function $u_0\in \mathcal{H}_{s,p}$ with  $\|u_0\|_{\mathcal{H}_{s,p}}>r>0$ such that \begin{align*}
\inf_{\|u\|_{\mathcal{H}_{s,p}}=r}E_{s,p} (u) >E_{s,p} (0)=0>E_{s,p} (u_0)
\end{align*}
for some positive constant $r$.
Secondly, this together with Lemma \ref{L3.4-2} implies that $E_{s,p} $ satisfies all assumptions of the
mountain-pass theorem. 
Define a set of paths
\begin{align*}
\Gamma=\left \{\gamma \in C([0,1], \mathcal{H}_{s,p}): \gamma(0)=0, \,\gamma(1)=u_0\right \}.
\end{align*}
 The mountain-pass theorem yields a critical point $u\in\mathcal H_{s,p}$
such that $E_{s,p}'(u)=0$ and
$$
 E_{s,p}(u)=
 \inf_{\gamma\in\Gamma}\max_{t\in[0,1]}
 E_{s,p}(\gamma(t))>0.
$$
Hence $u\not\equiv0$ is a weak solution of the truncated equation, and Lemma \ref{L4.2} implies
that $u$ is strictly positive and satisfies the original equation \eqref{eS1}.
\end{proof}

\subsection{Ground state solutions}\label{sub3}

Throughout this subsection, we  assume that $f$ also satisfies
$(\mathrm{A}_5)$.
  Define the Nehari manifold as $\mathcal{N}_{s,p}=\{u \in \mathcal{H}_{s,p}\backslash\{0\}:
 \langle E_{s,p}' (u), u\rangle =0  \}$, namely
	\begin{align}\label{N-1}
		\mathcal{N}_{s,p}=\left\{u \in \mathcal{H}_{s,p}\backslash\{0\}:   \int_{V} (|\nabla^s u|^p+h|u|^p ) d \mu=\int_{V}f(x,u^+ )u^+  d\mu \right\}.
	\end{align}
It is obvious that   $u^+ \not\equiv 0$ and $\|u\|_{\mathcal{H}_{s,p}}>0$ for any $u\in\mathcal{N}_{s,p}$.
By   assumption $(\mathrm{A}_4)$,
we get  two sufficiently small constants $\epsilon\in(0,\lambda_p)$   and $\delta>0$ such that
\begin{align} \label{e24}
 f(x, y)  \leq  (\lambda_p-\epsilon) y^{p-1}, \quad  \forall  x\in V,\, y\in[0,\delta).
\end{align}

We claim that
\begin{align}\label{e20}
\|u^+\| _{L^\infty } \geq \delta,\quad \forall u\in 	\mathcal{N}_{s,p}.
\end{align}
Otherwise, there  exists $u\in\mathcal{N}_{s,p}$
such that  $\|u^+\|_{L^\infty}<\delta$.
The preceding estimate gives
$$
\begin{aligned}
\|u\|_{\mathcal{H}_{s,p}}^p
&=\int_V f(x,u^+)u^+\,d\mu\\
&\leq(\lambda_p-\epsilon)\|u^+\|_{L^p}^p\\
&\leq\frac{\lambda_p-\epsilon}{\lambda_p}
       \|u\|_{\mathcal{H}_{s,p}}^p,
\end{aligned}
$$
which is impossible because $u\not\equiv 0$.

Moreover, we define
	\begin{align}\label{c}
		e_{s,p}=\inf_{u\in\mathcal{N}_{s,p}}E_{s,p} (u).
	\end{align}
 To construct a ground state solution of \eqref{eS1-1}, we first show that
$e_{s,p}$ is attained on $\mathcal N_{s,p}$ and then prove
that the resulting minimizer is a critical point of $E_{s,p}$.

\begin{lemma}\label{L4}
For any $  u \in \mathcal{H}_{s,p}$  with $u^{+} \not\equiv 0$,
 there exists a unique $t_0 >0$ such that
$t_0 u \in \mathcal{N}_{s,p}$ and
$E_{s,p} (t_0 u)=\underset{t>0}\max \,E_{s,p} (t u). $
  Moreover, if $u\in \mathcal{N}_{s,p}$, then $E_{s,p} ( u)=\underset{t>0}\max \,E_{s,p} (t u)$.
\end{lemma}

\begin{proof}
Fix $u\in\mathcal{H}_{s,p}$ with $u^+\not\equiv0$.
For any $t>0$,
let
\begin{align*}
  \phi (t)=\|u\|_{\mathcal{H}_{s,p}}^p-\int_{V}  u^+\frac{f(x, tu^+ ) }{t^{p-1}}    d\mu.
  \end{align*}
From Lemma \ref{Le4-1},   $\mathcal{H}_{s,p}$ is continuously embedded in $L^p(V)\cap L^\infty(V)$, and then   $u\in L^p(V)\cap L^\infty(V)$.
Assumptions $(\mathrm{A}_1)$--$(\mathrm{A}_4)$
and the dominated convergence theorem imply that
$\phi$ is well defined and continuous on $(0,\infty)$.
Moreover,  $(\mathrm{A}_5)$ and $u^+\not\equiv0$
imply that $\phi$ is strictly decreasing.

 Next, let $\alpha>p$ be as in $(\mathrm{A}_3)$.
By $(\mathrm{A}_3)$ and \eqref{eq:F-scaling}, for any $t\ge1$, we obtain
\begin{align*}
u^+\frac{f(x,tu^+)}{t^{p-1}}
\geq \frac{\alpha}{t^p}F(x,tu^+)
\geq \alpha t^{\alpha-p}F(x,u^+),
 \end{align*}
  and then
\begin{align}\label{e16}
\phi(t)
\leq \|u\|_{\mathcal{H}_{s,p}}^p
-\alpha t^{\alpha-p}\int_V F(x,u^+)\,d\mu .
\end{align}
Moreover, $(\mathrm{A}_3)$ and $u^+\not\equiv0$ give
$$
0<\int_V F(x,u^+)\,d\mu
\leq 
 \frac{\|u\|_{\mathcal{H}_{s,p}}^p-\phi(1)}{\alpha}
<\infty.
$$
This together with $\alpha>p$ and \eqref{e16} leads to
$\phi  (t)\rightarrow-\infty$ as $t \rightarrow\infty$.

 Finally, for all sufficiently small $t>0$, we have
$0\le t u^+(x)<\delta$ for every $x\in V$, where $\delta$ is chosen in \eqref{e24}.
Consequently,
$$
\phi(t)
\ge \|u\|_{\mathcal H_{s,p}}^p
 -(\lambda_p-\epsilon)\|u^+\|_{L^p}^p
\ge \frac{\epsilon}{\lambda_p}\|u\|_{\mathcal H_{s,p}}^p>0.
$$
By continuity and strict monotonicity, there exists a unique $t_0>0$
such that $\phi(t_0)=0$.
 Since
\begin{align}\label{6}
 \langle E_{s,p}'(tu),u\rangle
 =t^{p-1}\phi(t),
\end{align}
we have  $t_0 u \in \mathcal{N}_{s,p}$   and   $E_{s,p} (t_0 u)=\max_{t>0}E_{s,p} (t u)$.
Moreover,  if $u\in\mathcal{N}_{s,p}$, then $u^+\not\equiv0$ and
 $\phi(1)=0$. The uniqueness of $t_0$ gives $t_0=1$.
This completes the proof.
\end{proof}

Secondly, we prove that the infimum $e_{s,p}$ is attained on $\mathcal N_{s,p}$.

\begin{lemma}\label{L3}
There exists a function $u_0 \in \mathcal{N}_{s,p}$ such that $E_{s,p} (u_0)=e_{s,p}>0.$
\end{lemma}

\begin{proof}
By \eqref{c}, there exists
$\{u_n\}\subseteq\mathcal N_{s,p}$ such that
$E_{s,p}(u_n)\to e_{s,p}$    as ${n\rightarrow\infty}$.
Using $(\mathrm{A}_3)$ and \eqref{N-1}, we obtain
\begin{align*}
E_{s,p}(u_n)
\geq
\left(\frac1p-\frac1\alpha\right)
\|u_n\|_{\mathcal H_{s,p}}^p.
\end{align*}
Since $\alpha>p$, the sequence $\{u_n\}$ is bounded
in $\mathcal H_{s,p}(V)$.
By Lemma~\ref{Le4-1},   there exists
$u_0\in\mathcal H_{s,p}(V)$ such that
$u_n\rightharpoonup u_0$ weakly in $\mathcal H_{s,p}(V)$
and $u_n\to u_0$ strongly in $L^q(V)$
for every $q\in[p,\infty]$.
Moreover,  it follows from \eqref{e20} that  $\|u_n^+\| _{L^\infty}  \geq \delta$, where the sufficiently small
constant    $\delta>0$   is given by \eqref{e24}.
Since
$$
\|u_n^+-u_0^+\|_{L^\infty}
\leq\|u_n-u_0\|_{L^\infty}\rightarrow0,
$$
we obtain $\|u_0^+\|_{L^\infty}\geq\delta>0$,
and thus $u_0^+\not\equiv0$.

Choose $M>\delta$ such that
$\|u_n\|_{L^\infty}\leq M$ for all $n$ and
$\|u_0\|_{L^\infty}\leq M$.
By $(\mathrm{A}_1)$--$(\mathrm{A}_4)$, there exists
$K_M>0$ such that
\begin{align}\label{2}
|f(x,y)|\leq K_My^{p-1},
\quad \forall x\in V,\, y\in [ 0, M].
\end{align}
Consequently, we have
\begin{align*}
|F(x,u_n^+)-F(x,u_0^+)|
\leq
K_M\bigl(|u_n|^{p-1}+|u_0|^{p-1}\bigr)|u_n-u_0|.
\end{align*}
By H\"older's inequality and the strong $L^p$ convergence,
\begin{align*}
\int_V|F(x,u_n^+)-F(x,u_0^+)|\,d\mu\leq
K_M\left(
\|u_n\|_{L^p}^{p-1}+\|u_0\|_{L^p}^{p-1}
\right)\|u_n-u_0\|_{L^p}
\to0.
\end{align*}
In particular, we obtain
\begin{align}\label{20}
\lim_{n\to\infty}\int_VF(x,u_n^+)\,d\mu
=
\int_VF(x,u_0^+)\,d\mu.
\end{align}

 Set $p'=p/(p-1)$.
On every finite subset of $V$, the pointwise convergence
of $u_n$ and $(\mathrm{A}_1)$ imply $L^{p'}$ convergence
of the nonlinear terms.
The estimate \eqref{2} and the strong $L^p$ convergence
provide uniform control of their $L^{p'}$ tails.
Combining these two observations, we obtain
\[
f(\cdot,u_n^+)\to f(\cdot,u_0^+)
\quad\text{strongly in }L^{p'}(V).
\]
H\"older's inequality gives
\begin{align*}
& \left|
\int_V f(x,u_n^+)u_n^+\,d\mu
-\int_V f(x,u_0^+)u_0^+\,d\mu
\right|\\
\leq&
\left\|f(\cdot,u_n^+)-f(\cdot,u_0^+)\right\|_{L^{p'}}
\left\|u_n^+\right\|_{L^p} +
\left\|f(\cdot,u_0^+)\right\|_{L^{p'}}
\left\|u_n^+-u_0^+\right\|_{L^p}.
\end{align*}
Since $\{u_n^+\}$ is bounded in $L^p(V)$ and
\begin{align*}
\|u_n^+-u_0^+\|_{L^p}
\leq\|u_n-u_0\|_{L^p}\to0,
\end{align*}
we conclude that
\begin{align}\label{7}
\int_V f(x,u_n^+)u_n^+\,d\mu
\to
\int_V f(x,u_0^+)u_0^+\,d\mu.
\end{align}

 It follows from    $u_n \in \mathcal{N}_{s,p}$   and \eqref{7} that
\begin{align}\label{17}
\|u_0\|_{\mathcal H_{s,p}}^p
&\leq
\liminf_{n\to\infty}\|u_n\|_{\mathcal H_{s,p}}^p
\nonumber\\
&=
\lim_{n\to\infty}
\int_Vf(x,u_n^+)u_n^+\,d\mu
\nonumber\\
&=
\int_Vf(x,u_0^+)u_0^+\,d\mu.
\end{align}
By weak lower semicontinuity and \eqref{20}, we have
\begin{align}\label{21}
E_{s,p}(u_0)
&\leq
\liminf_{n\to\infty}
\left(
\frac1p\|u_n\|_{\mathcal H_{s,p}}^p
-\int_VF(x,u_n^+)\,d\mu
\right)
\nonumber\\
&=
\lim_{n\to\infty}E_{s,p}(u_n)
=e_{s,p}.
\end{align}

We claim that $u_0\in\mathcal N_{s,p}$.
Otherwise, since $u_0^+\not\equiv0$,
\eqref{17} and the definition of $\mathcal N_{s,p}$ imply
\begin{align}\label{N-strict-inequality}
\|u_0\|_{\mathcal H_{s,p}}^p
<
\int_Vf(x,u_0^+)u_0^+\,d\mu.
\end{align}
By Lemma~\ref{L4}, there exists a unique $t_0>0$
such that $t_0u_0\in\mathcal N_{s,p}$.
Thus $$E_{s,p}(t_0u_0)\geq e_{s,p}.$$
Recalling $u_n \in \mathcal{N}_{s,p}$ and  Lemma \ref{L4}, we have \begin{align*}
E_{s,p} (u_n)=\max _{t >0} E_{s,p} (tu_n).
\end{align*}
Since $t_0>0$ is fixed, $t_0u_n\to t_0u_0$ strongly
in $L^p(V)$, and these functions have a common
$L^\infty$ bound.
The same argument as for \eqref{20} therefore yields
\begin{align}\label{F-scaled-convergence}
\int_VF(x,t_0u_n^+)\,d\mu
\to
\int_VF(x,t_0u_0^+)\,d\mu.
\end{align}
Using \eqref{7}, \eqref{N-strict-inequality},
\eqref{F-scaled-convergence}, the Nehari identity for $u_n$,
and Lemma~\ref{L4}, we obtain
\begin{align*}
e_{s,p}
&\leq E_{s,p}(t_0u_0)\\
&=
\frac{t_0^p}{p}\|u_0\|_{\mathcal H_{s,p}}^p
-\int_VF(x,t_0u_0^+)\,d\mu\\
&<
\frac{t_0^p}{p}\int_Vf(x,u_0^+)u_0^+\,d\mu
-\int_VF(x,t_0u_0^+)\,d\mu\\
&=
\lim_{n\to\infty}E_{s,p}(t_0u_n)\\
&\leq
\lim_{n\to\infty}E_{s,p}(u_n)
=e_{s,p}.
\end{align*}
This contradiction shows that $u_0\in\mathcal N_{s,p}$.

Since $u_0\in\mathcal{N}_{s,p}$, \eqref{c} and  \eqref{21} give
$E_{s,p}(u_0)=e_{s,p}$.
Now $(\mathrm{A}_3)$ and the Nehari identity yield
$$
e_{s,p}=E_{s,p}(u_0)
\geq\left(\frac1p-\frac1\alpha\right)
       \|u_0\|_{\mathcal{H}_{s,p}}^p
>0.
$$
This completes the proof.
\end{proof}

By Lemmas~\ref{L4} and \ref{L3}, there exists
$u_0\in\mathcal N_{s,p}$ such that
\begin{align}\label{22-s}
\max_{t>0}E_{s,p} (t u_0)=E_{s,p} (u_0)=e_{s,p} >0.
\end{align}
We next show that $u_0$ is a critical point of $E_{s,p}$.

\begin{lemma}\label{L2}
We have $E_{s,p}' (u_0)=0$.
\end{lemma}

\begin{proof}
Motivated by \cite{S-Y-Z-2}, we adapt the variational
argument in \cite[Lemma~3.5]{Adimurthi}
to prove this lemma.
Suppose that $E_{s,p}'(u_0)\neq0$.
Choose $v\in\mathcal H_{s,p}(V)\setminus\{0\}$ such that
$\langle E_{s,p}'(u_0),v\rangle<0$, and define
\begin{align*}
\phi (t, m)=t u_0+m v, \quad \forall t ,    m \in \mathbb{R}.
\end{align*}
Since $E_{s,p}\in C^1(\mathcal H_{s,p}(V);\mathbb R)$,
we have
\begin{align*}
\frac{\partial}{\partial m}E_{s,p}(\phi(t,m))
=
\langle E_{s,p}'(\phi(t,m)),v\rangle
\to
\langle E_{s,p}'(u_0),v\rangle<0
\end{align*}
as $(t,m)\to(1,0)$.
Hence there exist $\epsilon_1\in(0,1)$ and
$\epsilon_2>0$ such that
\begin{align*}
\frac{\partial}{\partial m}E_{s,p}(\phi(t,m))
=
\langle E_{s,p}'(\phi(t,m)),v\rangle<0
\end{align*}
whenever $t\in[1-\epsilon_1,1+\epsilon_1]$ and
$m\in[-\epsilon_2,\epsilon_2]$.
Consequently,  \eqref{22-s} gives
\begin{align}\label{23}
E_{s,p}(\phi(t,m))
<
E_{s,p}(\phi(t,0))
=
E_{s,p}(tu_0)
\leq e_{s,p}
\end{align}
for every
$t\in[1-\epsilon_1,1+\epsilon_1]$ and
$m\in(0,\epsilon_2]$.

 For $(t,m)\in\mathbb R^2$, define
\begin{align}\label{N-perturbation-constraint}
\varphi(t,m)
 =&
\langle E_{s,p}'(\phi(t,m)),\phi(t,m)\rangle
\nonumber\\
 =&
\left\|\phi(t,m)\right\|_{\mathcal H_{s,p}}^p
-\int_V
f\bigl(x,\phi(t,m)^+\bigr)\phi(t,m)^+\,d\mu.
\end{align}
Since $E_{s,p}\in C^1(\mathcal H_{s,p}(V);\mathbb R)$,
the function $\varphi$ is continuous on $\mathbb R^2$.
Applying \eqref{6} to $u_0\in\mathcal N_{s,p}$
and using the strict monotonicity established in the proof
of Lemma~\ref{L4}, we obtain
\begin{align*}
\varphi(1+\epsilon_1,0)
<
0=\varphi(1,0)
<
\varphi(1-\epsilon_1,0).
\end{align*}
From continuity, there exists
$\epsilon_3\in(0,\epsilon_2)$ such that
\begin{align}\label{N-endpoint-signs}
 \varphi (1+\epsilon_1 ,m)  < 0< { \varphi (1-\epsilon_1 ,m)}, \quad \forall m \in\left(0, \epsilon_3\right].
\end{align}

 By decreasing $\epsilon_3$ if necessary, we may further
assume that
$$
\epsilon_3\|v\|_{\mathcal H_{s,p}}
<
\frac{1-\epsilon_1}{2}
\|u_0\|_{\mathcal H_{s,p}}.
$$
For every $t\in[1-\epsilon_1,1+\epsilon_1]$
and $|m|\leq\epsilon_3$, the reverse triangle inequality
then gives
\begin{align}\label{10}
 \|\phi(t,m)\|_{\mathcal H_{s,p}}
& =\|tu_0+mv\|_{\mathcal H_{s,p}} \nonumber\\
& \geq
t\|u_0\|_{\mathcal H_{s,p}}
-|m|\|v\|_{\mathcal H_{s,p}} \nonumber\\
& \geq
(1-\epsilon_1)\|u_0\|_{\mathcal H_{s,p}}
-\epsilon_3\|v\|_{\mathcal H_{s,p}} \nonumber\\
& >
\frac{1-\epsilon_1}{2}\|u_0\|_{\mathcal H_{s,p}}
>0.
\end{align}

 Fix $m\in(0,\epsilon_3]$.
By the continuity of $t\mapsto\varphi(t,m)$ and
\eqref{N-endpoint-signs}, the intermediate value theorem
yields $t_m\in(1-\epsilon_1,1+\epsilon_1)$ such that
$\varphi(t_m,m)=0$.
By \eqref{10}, we have $\phi(t_m,m)\neq0$.
Thus \eqref{N-perturbation-constraint} implies
$\phi(t_m,m)\in\mathcal N_{s,p}$.
Combining this with \eqref{23}, we obtain the   contradiction
\begin{align*}
e_{s,p}
\leq E_{s,p}(\phi(t_m,m))
<e_{s,p}.
\end{align*}
This completes the proof.
\end{proof}

\begin{proof}[\textbf{Proof of Theorem  \ref{T4}}]
 By Lemmas~\ref{L4}--\ref{L2}, there exists $u_0\in\mathcal{N}_{s,p}$ such that
$$
E_{s,p}(u_0)=e_{s,p}
=\inf_{u\in\mathcal{N}_{s,p}}E_{s,p}(u),
\quad E'_{s,p}(u_0)=0.
$$
Since every nontrivial critical point of $E_{s,p}$ belongs to
$\mathcal N_{s,p}$, $u_0$ is a ground state solution of \eqref{eS1-1}.

 By $(\mathrm{A}_5)$, we have
$f(x,t)=f(x,t^+)$ and $F(x,t)=F(x,t^+)$ for all $x\in V$ and $t\in\mathbb{R}.$
Therefore, \eqref{eS1} and \eqref{eS1-1} coincide,
and their associated energy functionals are identical.
Consequently, $u_0$ is also a ground state solution of
\eqref{eS1}. Moreover, Lemma~\ref{L4.2} implies that
$u_0$ is strictly positive in $V$.
This completes the proof.
\end{proof}



\begin{thebibliography}{00}



\bibitem{A1}
Adams R A. Sobolev spaces.
Pure and Applied Mathematics, Vol. 65.
New York--London: Academic Press, 1975.

\bibitem{Adimurthi}
Adimurthi. Existence of positive solutions of the semilinear Dirichlet problem with critical growth for the $n$-Laplacian. Ann. Scuola Norm. Sup. Pisa Cl. Sci. (4), 1990, 17: 393--413.

\bibitem{A-R}
Ambrosetti A, Rabinowitz P H.
Dual variational methods in critical point theory and applications.
J. Funct. Anal., 1973, 14: 349--381.

\bibitem{30}
Baleanu D, Diethelm K, Scalas E, et al.
Fractional calculus: Models and numerical methods.
Series on Complexity, Nonlinearity and Chaos, Vol. 3.
Singapore: World Scientific, 2012.

\bibitem{29}
Bogdan K, Byczkowski T. Potential theory for the $\alpha$-stable Schr\"{o}dinger operator on bounded Lipschitz domains. Studia Math., 1999, 133: 53--92.

\bibitem{A12}
Brezis H, Mironescu P. Gagliardo--Nirenberg, composition and products in fractional Sobolev spaces. J. Evol. Equ., 2001, 1: 387--404.

\bibitem{CRSTV2015}
Ciaurri \'{O}, Lizama C, Roncal L, et al. On a connection between the discrete fractional Laplacian and superdiffusion. Appl. Math. Lett., 2015, 49: 119--125.

\bibitem{CRSTV2018}
Ciaurri \'{O}, Roncal L, Stinga P R, et al.
Nonlocal discrete diffusion equations and the fractional discrete
Laplacian, regularity and applications.
Adv. Math., 2018, 330: 688--738.

\bibitem{De-C}
De Coster C, Dovetta S, Galant D, et al. On the notion of ground state for nonlinear Schr\"{o}dinger equations on metric graphs. Calc. Var. Partial Differential Equations, 2023, 62: Paper No. 159, 28 pp.

\bibitem{del-2021}
del Teso F, G\'{o}mez-Castro D, V\'{a}zquez J L.
Three representations of the fractional $p$-Laplacian:
semigroup, extension and Balakrishnan formulas.
Fract. Calc. Appl. Anal., 2021, 24: 966--1002.

\bibitem{A}
Di Nezza E, Palatucci G, Valdinoci E. Hitchhiker's guide to the fractional Sobolev spaces. Bull. Sci. Math., 2012, 136: 521--573.

\bibitem{Ge}
Ge H. Kazdan--Warner equation on graph in the negative case. J. Math. Anal. Appl., 2017, 453: 1022--1027.

\bibitem{28}
Getoor R K.
First passage times for symmetric stable processes in space.
Trans. Amer. Math. Soc., 1961, 101: 75--90.

\bibitem{A-Y-Y-2}
Grigor'yan A, Lin Y, Yang Y.
Kazdan--Warner equation on graph.
Calc. Var. Partial Differential Equations,
2016, 55: Paper No. 92, 13 pp.

\bibitem{A-Y-Y-1}
Grigor'yan A, Lin Y, Yang Y. Yamabe type equations on graphs. J. Differential Equations, 2016, 261: 4924--4943.

\bibitem{A-Y-Y-3}
Grigor'yan A, Lin Y, Yang Y. Existence of positive solutions to some nonlinear equations on locally finite graphs. Sci. China Math., 2017, 60: 1311--1324.

\bibitem{Haeseler-Keller-Lenz-Wojciechowski}
Haeseler S, Keller M, Lenz D, et al. Laplacians on infinite graphs: Dirichlet and Neumann boundary conditions. J. Spectr. Theory, 2012, 2: 397--432.

\bibitem{Han-Shao}Han X, Shao M. $p$-Laplacian equations on locally finite graphs. Acta Math. Sin. (Engl. Ser.), 2021, 37: 1645--1678.

\bibitem{Han-Shao-Zhao2020}
Han X, Shao M, Zhao L.
Existence and convergence of solutions for nonlinear
biharmonic equations on graphs.
J. Differential Equations, 2020, 268: 3936--3961.

\bibitem{HS}
Hou S, Sun J. Existence of solutions to Chern--Simons--Higgs equations on graphs. Calc. Var. Partial Differential Equations, 2022, 61: Paper No. 139, 13 pp.

\bibitem{Hua-Lin2017}
Hua B, Lin Y.
Stochastic completeness for graphs with curvature
dimension conditions.
Adv. Math., 2017, 306: 279--302.

\bibitem{Hua-Xu}
Hua B, Xu W. Existence of ground state solutions to some nonlinear Schr\"{o}dinger equations on lattice graphs. Calc. Var. Partial Differential Equations, 2023, 62: Paper No. 127, 17 pp.

\bibitem{Huang}
Huang X. On stochastic completeness of weighted graphs.
Ph.D. thesis. Bielefeld: Bielefeld University, 2011.

\bibitem{Keller-Lenz}
Keller M, Lenz D. Dirichlet forms and stochastic completeness of graphs and subgraphs. J. Reine Angew. Math., 2012, 666: 189--223.

\bibitem{Keller-1}
Keller M, Lenz D, Wojciechowski R K.
Graphs and discrete Dirichlet spaces.
Grundlehren der mathematischen Wissenschaften, Vol. 358.
Cham: Springer, 2021.

\bibitem{KN23}
Keller M, Nietschmann M. Optimal Hardy inequality for fractional Laplacians on the integers. Ann. Henri Poincar\'{e}, 2023, 24: 2729--2741.

\bibitem{Kwa15}
Kwa\'{s}nicki M. Ten equivalent definitions of the fractional Laplace operator. Fract. Calc. Appl. Anal., 2017, 20: 7--51.

\bibitem{Liu-Zhang}
Liu Y, Zhang M. Existence of solutions for nonlinear biharmonic Choquard equations on weighted lattice graphs. J. Math. Anal. Appl., 2024, 534: Paper No. 128079, 18 pp.

\bibitem{LM22}
Lizama C, Murillo-Arcila M. On a connection between the $N$-dimensional fractional Laplacian and $1$-D operators on lattices. J. Math. Anal. Appl., 2022, 511: Paper No. 126051, 12 pp.

\bibitem{LR2018}
Lizama C, Roncal L. H\"{o}lder--Lebesgue regularity and almost periodicity for semidiscrete equations with a fractional Laplacian. Discrete Contin. Dyn. Syst., 2018, 38: 1365--1403.

\bibitem{Peral}
Peral I. Multiplicity of solutions for the $p$-Laplacian.
Lecture notes, Second School of Nonlinear Functional Analysis
and Applications to Differential Equations.
Trieste: ICTP, 1997.

\bibitem{A78}
Runst T, Sickel W. Sobolev spaces of fractional order, Nemytskij operators, and nonlinear partial differential equations. Berlin: De Gruyter, 1996.

\bibitem{S-Y-Z-1}
Shao M, Yang Y, Zhao L. Sobolev spaces on locally finite graphs. Proc. Amer. Math. Soc., 2025, 153: 693--708.

\bibitem{S-Y-Z-2}
Shao M, Yang Y, Zhao L. Existence and convergence of solutions to $p$-Laplace equations on locally finite graphs. J. Korean Math. Soc., 2026, 63: 821--848.

\bibitem{S-w}
Sun L, Wang L. Brouwer degree for Kazdan--Warner equations on a connected finite graph. Adv. Math., 2022, 404: Paper No. 108422, 29 pp.

\bibitem{Tian-Zhao2025}
Tian Y, Zhao L.
A note for $W^{1,p}(V)$ and $W^{1,p}_0(V)$
on a locally finite graph.
Acta Math. Sci. Ser. B (Engl. Ed.),
2025, 45: 2135--2141.

\bibitem{31}
Valdinoci E. From the long jump random walk to the fractional Laplacian. Bol. Soc. Esp. Mat. Apl. SeMA, 2009, 49: 33--44.

\bibitem{Wang-frac}
Wang J. Eigenvalue estimates for the fractional Laplacian on lattice subgraphs. arXiv:2303.15766, 2023.

\bibitem{Yu2026}
Yu P. Existence for nonlinear fractional $p$-Laplacian
equations on finite graphs.
J. Korean Math. Soc., 2026, 63: 279--293.

\bibitem{Z-L-Y1}
Zhang M, Lin Y, Yang Y. Fractional Laplace operator on finite graphs.
Rev. Mat. Complut., 2026, 39: 721--753.

\bibitem{Z-L-Y2}
Zhang M, Lin Y, Yang Y. Fractional Laplace operator and related Schr\"{o}dinger equations on locally finite graphs. Calc. Var. Partial Differential Equations, 2025, 64: Paper No. 227, 27 pp.

\bibitem{Zhang-Zhao}Zhang N, Zhao L. Convergence of ground state solutions for nonlinear Schr\"{o}dinger equations on graphs. Sci. China Math., 2018, 61: 1481--1494.
\end{thebibliography}
\end{document}